\documentclass[aos,citesort,dvips]{arximspdf}
\usepackage{dcolumn}
\usepackage{graphicx}

\doi{10.1214/10-AOS793}
\volume{38}
\issue{5}
\pubyear{2010}
\firstpage{2620}
\lastpage{2651}

\makeatletter

\newcolumntype{d}[1]{D{.}{.}{#1}}

\newtheorem{theorem}{Theorem}
\newproclaim{remark}{Remark}
\newtheorem{lemma}{Lemma}
\newproclaim{assumption}{Assumption}

\makeatother

\begin{document}
\begin{frontmatter}

\title{Sparse recovery under matrix uncertainty}
\runtitle{Sparse recovery under matrix uncertainty}

\begin{aug}
\author[A]{\fnms{Mathieu} \snm{Rosenbaum}\ead[label=e1]{mathieu.rosenbaum@polytechnique.edu}} and
\author[B]{\fnms{Alexandre B.} \snm{Tsybakov}\corref{}\thanksref{T3}\ead[label=e2]{alexandre.tsybakov@upmc.fr}}
\runauthor{M. Rosenbaum and A. B. Tsybakov}
\affiliation{CMAP-\'{E}cole Polytechnique Paris, CREST and
LPMA-University of Paris~6}
\address[A]{CMAP-\'{E}cole Polytechnique Paris\\
UMR CNRS 7641\\
91128 Palaiseau Cedex\\
France\\
\printead{e1}}
\address[B]{Laboratoire de Statistique, CREST\\
3, av. Pierre Larousse\\
92240 Malakoff\\
France \\
and\\
LPMA (UMR CNRS 7599)\\
Universit\'{e} Paris 6\\
4, Place Jussieu\\
75252 Paris, Cedex 05\\
France\\
\printead{e2}}
\end{aug}

\thankstext{T3}{Supported in part by the Grant ANR-06-BLAN-0194 and by
the PASCAL
Network of Excellence.}

\received{\smonth{12} \syear{2008}}
\revised{\smonth{10} \syear{2009}}

%
\begin{abstract}
We consider the model
\begin{eqnarray*}
y&=& X\theta^* + \xi,\\
Z&=&X+\Xi,
\end{eqnarray*}
where the random vector $y\in\mathbb{R}^n$ and the random $n\times p$
matrix $Z$ are observed, the $n\times p$ matrix $X$ is unknown,
$\Xi$ is an $n\times p$ random noise matrix, $\xi\in\mathbb{R}^n$
is a
noise independent of $\Xi$, and $\theta^*$ is a vector of unknown
parameters to be estimated. The matrix uncertainty is in the fact
that $X$ is observed with additive error. For dimensions $p$ that
can be much larger than the sample size $n$, we consider the
estimation of sparse vectors $\theta^*$. Under matrix uncertainty,
the Lasso and Dantzig selector turn out to be extremely unstable
in recovering the sparsity pattern (i.e., of the set of nonzero
components of $\theta^*$), even if the noise level is very small.
We suggest new estimators called \textit{matrix uncertainty
selectors} (or, shortly, the \textit{MU-selectors}) which are close
to $\theta^*$ in different norms and in the prediction risk if the
restricted eigenvalue assumption on $X$ is satisfied. We also show
that under somewhat stronger assumptions, these estimators recover
correctly the sparsity pattern.
\end{abstract}

%
\begin{keyword}[class=AMS]
\kwd[Primary ]{62J05}
\kwd[; secondary ]{62F12}.
\end{keyword}
\begin{keyword}
\kwd{Sparsity}
\kwd{MU-selector}
\kwd{matrix uncertainty}
\kwd{errors-in-variables model}
\kwd{measurement error}
\kwd{sign consistency}
\kwd{oracle inequalities}
\kwd{restricted eigenvalue assumption}
\kwd{missing data}
\kwd{portfolio selection}
\kwd{portfolio replication}.
\end{keyword}

\end{frontmatter}

\section{Introduction}\label{intro}

We consider the model
%
\begin{eqnarray}\label{0}
y&=& X\theta^* + \xi,\\
\label{00}
Z&=&X+\Xi,
\end{eqnarray}
where the random vector $y\in\mathbb{R}^n$ and the random $n\times p$
matrix $Z$ are observed, the $n\times p$ matrix $X$ is unknown,
$\Xi$ is an $n\times p$ random noise matrix, $\xi\in\mathbb{R}^n$
is a
noise independent of $\Xi$, and
$\theta^*=(\theta^*_1,\ldots,\theta^*_p)$ is a vector of unknown
parameters to be estimated.

We will typically assume that $\theta^*$ is $s$-sparse, that is, that
it has only $s$ nonzero components, where $1\le s \le p$ is some
integer. The dimension $p$ can be much larger than the sample size
$n$, but we will typically have in mind the situation where the
effective dimension $s$ is much smaller than $p$ and $n$. We will
also assume that the elements of $\Xi$ are small. In this setting
we will suggest estimators
$\hat\theta=(\hat\theta_1,\ldots,\hat\theta_p)$ that under some
assumptions recover $\theta^*$ with high accuracy in different
norms, as well as under the prediction risk. We will also show
that, under somewhat stronger assumptions, these estimators recover
correctly the sparsity pattern, that is, the set of nonzero
components of $\theta^*$. Our results follow the spirit of the now
extensive literature on sparsity with $\ell_1$-minimization (see,
e.g.,
\cite{brt,btw07a,btw07b,c08,crt,ct1,ct,det,l,l2,kol,kolsf,mb,vdg08,zh,zhT,zy,zou}).
The main difference is in the presence of matrix uncertainty. The
matrix $X$ is not known and is observed with error. This leads us
to new estimators, called \textit{matrix uncertainty selectors} (or,
shortly, the MU-selectors), which are different from the Lasso
and Dantzig selector (or their modifications) studied in those
papers.

In what follows, without loss of generality, we mainly assume that
$\xi$ and $\Xi$ are deterministic and satisfy the assumptions
%
\begin{eqnarray}
\label{1}
\biggl|\frac1{n}Z^T\xi\biggr|_\infty&\le& \varepsilon,\\
\label{2}
|\Xi|_\infty&\le& \delta
\end{eqnarray}
for some $\varepsilon\ge0, \delta\ge0$ [a modification of (\ref
{2}) is also
used in some cases]. Here $ |\cdot|_\infty$ stands for
the maximum of components norm. If $\xi$ and $\Xi$ are random,
conditions (\ref{1}) and (\ref{2}) can be guaranteed with a
probability close to 1 under natural assumptions that we discuss
below; we also indicate the corresponding values of $\varepsilon$
and~$\delta$.
So, the results that we prove for deterministic $\xi$ and $\Xi$
are extended in a trivial way to random $\xi$ and $\Xi$ satisfying
these assumptions. The difference is only in the fact that the
results hold on the random event of high probability where
(\ref{1}) and (\ref{2}) are satisfied. The setting with random $X$
is covered in a similar way. We only need to consider random $X$
for which the restricted eigenvalue (RE) assumption or the
Coherence assumption (see below) hold with high probability.
Examples of such random $X$ are discussed in the literature
\cite{ct,kolsf}.

We introduce two versions of MU-selectors. The first one is
designed for the case $\xi=0$, that is, for the problem of solving a
large system of linear equations with deterministic or random noise
in the matrix. This MU-selector is defined as a solution of the
minimization problem
\[
\min\{ |\theta|_1\dvtx \theta\in\Theta, |y-Z\theta|_\infty\le
\delta
|\theta|_1 \},
\]
where $\Theta\subseteq\mathbb{R}^p$ is a given set characterizing the
prior knowledge about $\theta$. Here and below $|x|_q$, $q\ge1$,
denotes the $\ell_q$-norm of $x\in\mathbb{R}^d$ whatever is $d\ge1$.

The second version of the MU-selector is defined as a solution of the
minimization problem
%
\begin{equation}\label{mu_sel}
\min\biggl\{ |\theta|_1\dvtx \theta\in\Theta,
\biggl|\frac1{n}Z^T(y-Z\theta) \biggr|_\infty\le\lambda|\theta|_1
+\varepsilon\biggr\},
\end{equation}
where $\lambda\ge0$ is a factor responsible for matrix
uncertainty. If (\ref{2}) is assumed, we choose $\lambda$
depending on $\delta$ so that $\lambda>0$ for $\delta>0$ and
$\lambda=0$
for $\delta=0$. Note that if $\Theta=\mathbb{R}^p$ and there is no matrix
uncertainty, that is, $\delta=0$, this second MU-selector becomes the
Dantzig selector of \cite{ct} based on the data $(y,Z)$. We mainly
discuss the choice $\lambda(\delta)=(1+\delta)\delta$, which
corresponds to a
noise satisfying (\ref{2}). This can be also used for $\xi=0$ by
setting $\varepsilon=0$ in the definition. Nevertheless, for $\xi=0$ we
consider directly the first version of the MU-selector because it is
simpler and achieves better error bounds than for
$\lambda=(1+\delta)\delta$.

Note that using in model (\ref{0}) and (\ref{00}) the Lasso or Dantzig
selector with $Z$ instead of the true $X$ typically leads to
satisfactory results for the prediction loss when the noise $\Xi$ is
small enough. However, these methods are less efficient in
estimation of $\theta^*$ and they are especially unstable in
selection of the sparsity pattern (cf. Section~\ref{appli}). In
particular, they become quite sensitive to the values of $\theta^*$.
This is explained by the fact that the true $\theta^*$ is no longer
guaranteed to stay, with a probability close to 1, in the feasible
set of the Dantzig selector (which is also the set containing all
the Lasso solutions).

The second MU-selector differs from the Dantzig selector based on
the data $(y,Z)$ in that we ``penalize more'' by enlarging the
feasible band for $ |\frac1{n}Z^T(y-Z\theta) |_\infty$. Indeed,
setting $\Theta=\mathbb{R}^p$, a Lasso type analog of this MU-selector
can be defined as a solution of the convex minimization problem
\[
\min_{\theta\in\mathbb{R}^p} \biggl\{\frac1{n}|y-Z\theta|_2^2 +
\lambda_1|\theta|_1 +\lambda_2 |\theta|_1^2 \biggr\}
\]
with some $\lambda_1,\lambda_2>0$. To appreciate why there is a
similarity, note that for $\theta$ to achieve the minimum of such a
convex criterion, it is necessary and sufficient to have
%
\begin{eqnarray}\label{las1}
\biggl(\frac{1}{n}Z^T(y-Z\theta) \biggr)_j &=&
\frac{\lambda_1}{2}+\lambda_2 |\theta|_1 \operatorname
{sign}(\theta_j)\qquad
\mbox{if $\theta_j\neq0$},\nonumber\\[-8pt]\\[-8pt]
\biggl| \biggl(\frac{1}{n}{Z}^T(y-Z\theta) \biggr)_j \biggr| &\le&
\frac{\lambda_1}{2}+\lambda_2 |\theta|_1 \qquad\mbox{if $\theta_j=
0$},\nonumber
\end{eqnarray}
where the index $j$ designates the $j$th component of the
corresponding vector and $\operatorname{sign}(\theta_j)$ is the sign of
$\theta_j$.
Therefore, the set of possible solutions is ``tightly'' contained in
$ \{ \theta\in\mathbb{R}^p\dvtx |\frac1{n}Z^T(y-Z\theta) |_\infty
\le
\lambda_2|\theta|_1 +\lambda_1/2 \}$, which is the feasible set of
the MU-selector (\ref{mu_sel}). The analogy is thus in the same
spirit as between the Lasso and the Dantzig selector.

The results of this paper can be viewed in several perspectives.
First, we can interpret them as a new approach to the inference in
errors-in-variables models. The classical ways of treating these
models via some versions of least squares or of the method of
moments heavily depend on specific identifiability constraints that
are violated when $p\gg n$ \cite{full,ds}. Our approach is free of
such constraints and requires only a modest price, which is the
sparsity of the unknown vector of parameters. Also, on the
difference from the results in the conventional errors-in-variables
framework, we provide nonasymptotic bounds for the risks of the
estimators and guarantee the finite sample variable selection
property.

The second perspective is an extension of the theory of
$\ell_1$-based sparse recovery beyond the restricted
isometry/restricted eigenvalue conditions (cf. \cite{ct,brt}) that
are known to be too strong. We show that small perturbations of the
design matrix $X$ that bring these conditions to failure are in fact
not so dangerous, once the method of recovery is chosen in a proper
way (cf. Remark \ref{Remark4} below).

Finally, the third perspective is in developing simple and
efficient tools of sparse recovery for specific applications. We
mention here models with missing data, some financial models
(portfolio selection, portfolio replication) and inverse problems
with unknown operator. They are presented in the next section.

\section{Examples of application}\label{sec2}
Here we explain how several examples of application can be
described by model (\ref{0}) and (\ref{00}) with a sparse vector of
parameters~$\theta^*$.

1. \textit{Models with missing data.} Assume that the elements $Z_{ij}$
of matrix $Z$ satisfy
%
\begin{equation}\label{ex1}
Z_{ij} = X_{ij}\eta_{ij},
\end{equation}
where $X_{ij}$ are the elements of $X$ and $\eta_{ij}$ are i.i.d.
Bernoulli random variables taking value 1 with probability $1-\pi$
and 0 with probability $\pi$, $0<\pi<1$. The data $X_{ij}$ is
missing if $\eta_{ij}=0$, which happens with probability $\pi$. We
are mainly interested in the case of small $\pi$. In practice, it is
easy to estimate $\pi$ by the empirical probability of occurrences
of zeros in the sample of $Z_{ij}$, so it is realistic to assume
that $\pi$ is known. Note that we can rewrite (\ref{ex1}) in the
form
%
\begin{equation}\label{ex2}
Z_{ij}^\prime= X_{ij}+ \xi_{ij}^\prime,
\end{equation}
where $Z_{ij}^\prime=Z_{ij}/(1-\pi)$,
$\xi_{ij}^\prime=X_{ij}(\eta_{ij}-E(\eta_{ij}))/(1-\pi)$ and
$E(\cdot)$ denotes the expectation. Thus, we can reduce the model
with missing data (\ref{ex1}) to the form (\ref{00}) with matrix
$\Xi$ whose elements $\xi_{ij}^\prime$ are zero mean bounded random
variables. In this case assumption (\ref{2}) is fulfilled with $\delta$
which is not necessarily small, whereas the theoretical bounds
obtained below only make sense if $\delta$ is small enough.
Nevertheless, the variances of $\xi_{ij}^\prime$ are proportional to
$\pi$, and we will see in Section \ref{sec6} that, by modifying assumption
(\ref{2}), we obtain bounds for the MU-selector that are small if $\pi$ is
small.

2. \textit{Portfolio selection.} Brodie et al. \cite{daub} recently
argued that classical methods of portfolio selection are highly
unstable. As a remedy, they proposed an algorithm accounting for the
sparsity of portfolio weights and studied its numerical performance.
A different approach to sparse portfolio selection can be introduced
in our framework. Recall that in the traditional Markowitz portfolio
selection, the objective is to find a portfolio having minimal
variance return for a given expected return. This is stated as the
optimization problem
%
\begin{equation}\label{portf}
\min\biggl\{ \theta^T X \theta\dvtx \theta^T\mu=\beta, \sum_j\theta
_j=1, j=1,\ldots, p\biggr\},
\end{equation}
where $\theta=(\theta_1,\ldots,\theta_p)$ is the vector of
weights with $\theta_j$
representing the proportion of capital invested in the $j$th
asset, $X$ and $\mu$ are the covariance matrix and the vector of
expected returns of the different assets and $\beta$ is the
desired return of the portfolio.

Using Lagrange multipliers, problem (\ref{portf}) is reduced to
the solution of the linear equation $X\theta=a$ for some vector $a\in
\mathbb{R}^p$ depending on $\mu$ and $X$. However, neither the covariance
matrix $X$, nor the mean $\mu$ are available. Only their empirical
(noisy) versions are observed. Instead of $X$ we have a sample
covariance matrix $Z$, and instead of $a$ a vector of noisy
observations $y$. Thus, we are in the framework of model
(\ref{0}) and (\ref{00}) with $n=p$ (since $X$ is a square matrix).
Direct substitution of noisy values $Z$ and $y$ instead of $X$ and
$a$ leads, in general, to instability of the solution of the
linear equation because the dimension $p$ can be very high (often
500 assets or more) and $X$ can be either degenerate or with a
small minimal eigenvalue. The methods that we suggest below are
robust to the variations both of the matrix $X$ and of the
right-hand side $a$.

Another way of looking at sparse portfolio selection is to revise
the very problem (\ref{portf}). Note that minimizing $\theta^T X
\theta$,
where $X$ is the covariance matrix, is motivated by the fact that we
would like to get the portfolio with smallest ``dispersion.'' This
requirement looks quite natural as long as we remain in the world of
the classical second order statistics reasoning. The problem
(\ref{portf}) is similar in spirit to ``minimal variance
unbiased estimation,'' an old concept which is known to have serious
drawbacks. An alternative method would be to look for the
\textit{sparsest portfolio with a given daily return $\beta$} (we can also
consider weakly or monthly returns). The problem can be formalized
as follows. Let $X_{ij}$ be the return of the $j$th asset on day
$i$. The matrix of returns $X=(X_{ij})_{i,j}$ is typically observed
with measurement error. This error can be due to an incomplete
description of the assets. For example, the only available
quantities for the investor are often reduced to the open, high, low
and close prices, which leads, in particular, to asynchronous data
(especially when one deals with prices from markets belonging to
different time zones) and to an underestimation of the order book
effects. Indeed, to take into account the liquidity costs, an
investor should compute the returns, having in mind the order of
magnitude of the number of assets he may have in his portfolio.
However, such an accurate computation is only possible for the very
few investors having access to order book data. In the case of
nonstandard assets such as hedge funds, the measurement error can be
also due to uncertainty about the management costs, the rounding
approximations used and the way the returns are computed. Thus,
instead of $X$, we in fact observe some other matrix
$Z=(Z_{ij})_{i,j}$.

We are looking for the sparsest portfolio, that is, a portfolio that
solves the problem
%
\begin{equation}\label{portf0}
\min\biggl\{ |\theta|_0\dvtx X\theta=\mathbf{b}, \sum_j\theta_j=1, j=1,\ldots,
p\biggr\},
\end{equation}
where $\theta$ is the vector of the proportions of the wealth invested
in each asset, $|\theta|_0$~is the number of nonzero components of
$\theta$
and $\mathbf{b}\in\mathbb{R}^n$ is the vector with all the components equal
to $\beta$. It is important to note that the sparsest portfolio does
not necessarily contain a very small number of assets, in particular,
when $p$ is large. The minimization problem (\ref{portf0}) is
NP-hard, and the standard way to approximate it is to consider its
convex relaxation:
%
\begin{equation}\label{portf1}
\min\biggl\{ |\theta|_1\dvtx X\theta=\mathbf{b}, \sum_j\theta_j=1, j=1,\ldots,
p\biggr\}.
\end{equation}
This problem is already numerically solvable, but since $X$ is
observed with error, the solution can be unstable. We do not
necessarily recover the sparsest solution if we directly plug $Z$
instead of $X$ in (\ref{portf1}). A stable alternative that we
suggest below is given by solving
%
\begin{equation}\label{portf2}
\min\biggl\{ |\theta|_1\dvtx |\mathbf{b}-Z\theta|_{\infty}\le\delta|\theta|_1,
\sum_j\theta_j=1, j=1,\ldots, p\biggr\},
\end{equation}
where $\delta$ is an upper bound on the noise level in the matrix $X$.

3. \textit{Portfolio replication.} Replicating a portfolio, or at
least finding the type of assets in a portfolio, has become a very
challenging issue in the recent years, especially in the hedge
funds context. Indeed, replicating a hedge fund portfolio means
obtaining a Profit and Loss profile similar to those of the hedge
fund without investing in it (and so avoiding the usual drawbacks
of a hedge fund investment such as excessive fees, lack of
transparency, lack of liquidity, lack of capacity, etc.).
Replicating a portfolio can be done by retrieving the assets
belonging to the portfolio. This problem can be formalized through
model (\ref{0}) and (\ref{00}).

To fix ideas, suppose, for example, that we observe the daily
returns $y_i$, $i=1,\ldots, T$, of a portfolio. Moreover, assume
that the proportion of capital invested in each asset of the
portfolio is constant between day $1$ and day $T$. Then, we
theoretically have
\[
y_i=\sum_{j=1}^p\theta_j X_{ij},
\]
where $p$ is the total number of different assets in the portfolio,
$X_{ij}$ is the return of the $j$th asset belonging to the
portfolio on day $i$ and $\theta_j$ the proportion of capital
invested in it. As pointed out in the previous example, it is
natural to consider that the vector of the portfolio returns
$(y_i)_i$ and the matrix of the assets returns $X=(X_{ij})_{i,j}$
are observed with measurement error. Note that in this setup we can
also treat the case where $y_i$ and $X_{ij}$ are the absolute
returns (differences between the close price and the open price),
provided that we define $\theta_j$ as the (constant) quantity of the
$j$th asset in the portfolio.

To solve our problem, we could formally consider that any existing
asset or derivative can, in principle, belong to the portfolio. This
is of course not realistic. However, it is reasonable to assume that
the portfolio is rather sparse and that any asset used in the
portfolio has a behavior which is quite close to those of an asset
belonging to a restricted, given class of reference assets,
especially if this restricted class can still be very large. For
example, we will not put all the oil companies in the world in our
restricted class. Nevertheless, we suppose that if an oil company is
used in the portfolio and is not in the class, its returns profile
will look like the returns profile of another oil company which
belongs to the restricted class. Consequently, $Z$ will be made from
the daily returns of our reference assets. Indeed, for any asset in
the portfolio, we will assume that either it belongs to the assets
defining $Z$ or it ``resembles'' one of the assets defining $Z$.
Eventually, $Z$ can be seen as a noisy measurement of $X$ and, thus,
the problem is described by model (\ref{0}) and (\ref{00}). A numerical
illustration is given in Section \ref{appli}.

4. \textit{Inverse problems with unknown operator.} This setting has
been recently discussed by several authors
\cite{ek,cavh,cavr,hr,m07}. A typical problem is to recover an
unknown function $f$ that belongs to a Hilbert space $H$ based on a
noisy observation $Y$ of $Af$ where $A\dvtx H\to V$ is a linear operator
and $V$ is another Hilbert space. The observation $Y$ with values in
$V$ can be written as
%
\begin{equation}\label{inv}
Y= Af + \zeta,
\end{equation}
where $\zeta$ is a random variable (typically assumed Gaussian) with
values in $V$. Let $\{\phi_j\}_{j=1}^\infty$ and
$\{\psi_j\}_{j=1}^\infty$ be complete orthonormal bases in $H$ and
$V$, respectively. We can write $f = \sum_{j=1}^\infty\theta_j^*
\phi_j=\sum_{j=1}^p \theta_j^* \phi_j + r$ with some coefficients
$\theta_j^*$, where the integer $p$ is chosen very large, so that one
can consider the remainder term $r\in H$ as negligible. Therefore,
we can reduce the problem of estimating $f$ to that of recovering
the vector of coefficients $\theta^* = (\theta_1^*,\ldots,\theta_p^*).$
Introducing the scalar products $Y_i=(Y,\psi_i)$ and
$\xi_i=(\zeta,\psi_i)$, we obtain from (\ref{inv}) the following
sequence of real-valued observations:
\[
Y_i = \sum_{j=1}^p \theta_j^* (A\phi_j,\psi_i) + (Ar,\psi_i)+
\xi_i,\qquad
i=1,2,\ldots.
\]
If we consider here only the first $n$ observations, assume that
$(Ar,\psi_i)=0$ and define the matrix
$X= ((A\phi_j,\psi_i)_{i=1,\ldots,n, j=1,\ldots,p} )$, and
the vectors $y=(Y_1,\ldots,Y_n)$, $\xi=(\xi_1,\ldots,\xi_n)$, then we
get the linear model (\ref{0}). As discussed in
\cite{ek,cavh,cavr,hr,m07}, it is rather frequent in the
applications that the operator $A$ is not known, but its action on
any given function in $H$ can be observed with some noise. We
emphasize that in those papers the noise $\Xi$ is supposed to be
small. This is consistent with the strategy of performing many
repeated measurements of $(A\phi_j,\psi_i)$ for each pair $(i,j)$.
Thus, we have access to observations of the matrix
$X= ((A\phi_j,\psi_i)_{i=1,\ldots,n, j=1,\ldots,p} )$ with
some small noise and, therefore, we are in the framework of model
(\ref{0}) and (\ref{00}). The results obtained in
\cite{ek,cavh,cavr,hr,m07} consider the case $n=p$ and deal with
nondegenerate matrices $X$. This framework is not always
convenient, especially if $n$ and $p$ are very large. The approach
that we develop in this paper is more general in the sense that, for
example, if $n=p$ we can treat degenerate matrices $X$ that satisfy
some regularity assumptions. We also cover the case $p\gg n$, which
is a useful extension because by taking a large $p$ we can assure
that the residual $r$ is indeed negligible.


\section{Sparse solution of linear equations with noisy
matrix}\label{sec3}

In this section we consider the simplest case, $\xi=0$. Thus, we
solve the system of linear equations
\[
y=X\theta,
\]
where $X$ is an unknown matrix such that we can observe its noisy
values
\[
Z=X+\Xi,
\]
where $\Xi$ satisfies (\ref{2}).

Let $\Theta$ be a given convex subset of $\mathbb{R}^p$. We will
assume in
this section that there exists an $s$-sparse solution $\theta_s$ of
$y=X\theta$ such that $\theta_s\in\Theta$. Consider the estimator
$\hat\theta$ of $\theta_s$ defined as a solution of the following
minimization problem:
%
\begin{equation}\label{3}
\min\{ |\theta|_1\dvtx \theta\in\Theta, |y-Z\theta|_\infty\le
\delta
|\theta|_1 \}.
\end{equation}
%

Clearly, (\ref{3}) is a convex minimization problem. If $\Theta
=\mathbb{R}^p$
or if $\Theta$ is a linear subspace of $\mathbb{R}^p$ or a simplex
(the latter
case is interesting, e.g., in the context of portfolio
selection), then (\ref{3}) reduces to a linear programming problem.

Note that under assumption (\ref{2}) the feasible set of problem
(\ref{3})
\[
\Theta_1 = \{ \theta\in\Theta\dvtx |y-Z\theta|_\infty\le\delta
|\theta|_1 \}
\]
is nonempty. In fact, $\theta_s\in\Theta_1$ since
%
\begin{equation}\label{3a}
|y-Z\theta_s |_\infty= |\Xi\theta_s |_\infty
\le|\Xi|_\infty|\theta_s|_1 \le\delta|\theta_s|_1.
\end{equation}
Thus, there always exists a solution $\hat\theta$ of (\ref{3}). But it
is not necessarily unique. We will call solutions of (\ref{3}) the
matrix uncertainty selectors (or, shortly, MU-selectors).

To state our assumptions on $X$, we need some notation. For a vector
$\theta\in\mathbb{R}^p$ and a subset $J$ of $\{1,\ldots,p\}$, we
denote by $\theta_J$
the vector in $\mathbb{R}^p$ that has the same coordinates as $\theta
$ on the
set of indices $J$ and zero coordinates on its complement~$J^c$.

We will assume that the matrix $X$ satisfies the following
condition (\textit{restricted eigenvalue assumption} \cite{brt}):
\renewcommand{\theassumption}{RE($s$)}
\begin{assumption}\label{assumpREs}\hypertarget{assumpRE2s}
There exists $\kappa>0$
such that
\[
\min_{\Delta\ne0\dvtx |\Delta_{J^c}|_1\le
|\Delta_{J}|_1}\frac{|X\Delta|_2}{\sqrt{n}|\Delta_J|_2} \ge
\kappa
\]
for all subsets $J$ of $\{1,\ldots,p\}$ of cardinality $|J|\le s$.
\end{assumption}

A detailed discussion of this assumption can be found in
\cite{brt}. In particular, it is shown in \cite{brt} that the
restricted eigenvalue assumption is more general than several
other similar assumptions used in the sparsity literature
\cite{det,ct,zh}. One of such assumptions is the \textit{coherence
condition} \cite{det} that has the following form.
\renewcommand{\theassumption}{C}
\begin{assumption}\label{assumpC}
All the diagonal elements of the
matrix $\Psi= X^TX/n$ are equal to 1 and all its off-diagonal
elements $\Psi_{ij}, i\neq j$, satisfy the coherence condition:
${\max_{i\neq j}}|\Psi_{ij}|\le\rho$ with some $\rho<1$.
\end{assumption}

Note that Assumption \ref{assumpC} with $\rho<(3\alpha s)^{-1}$ implies Assumption
\ref{assumpREs} with $\kappa=\sqrt{1-1/\alpha}$ (cf. \cite{brt} or Lemma 2 in
\cite{l}).

We now state the main result of this section.
\begin{theorem}\label{t1}
Assume that there exists an $s$-sparse solution $\theta_s\in\Theta
$ of the
equation $y=X\theta$. Let (\ref{2}) hold.
Then for any solution $\hat\theta$
of (\ref{3}) we have the following inequalities:

\begin{longlist}
\item
%
\begin{equation}\label{6}
\frac{1}{n}|X(\hat\theta-\theta_s)|_2^2\le 4\delta^2|\hat
\theta
|_1^2.
\end{equation}
\item If Assumption \ref{assumpREs} holds, then
%
\begin{equation}\label{4}
|\hat\theta-\theta_s|_1\le \frac{4\sqrt{s}\delta}{\kappa
}|\hat\theta|_1.
\end{equation}
\item If Assumption \textup{RE($2s$)} holds, then
%
\begin{equation}\label{5}
|\hat\theta-\theta_s|_2\le
\frac{4\delta}{\kappa}|\hat\theta|_1.
\end{equation}
\item If Assumption \ref{assumpC} holds with $\rho<\frac1{3\alpha s}, \alpha
>1$, then
%
\begin{equation}\label{6a}
|\hat\theta-\theta_s|_\infty<
2 \biggl(1+\frac{2}{3\sqrt{s\alpha(\alpha-1)}} \biggr)\delta|\hat\theta|_1.
\end{equation}
\end{longlist}
\end{theorem}
\begin{pf}
Set $\Delta=\hat\theta-\theta_s$ and $J=J(\theta
_s)$, where
$J(\theta)$ denotes the set of nonzero coordinates of $\theta$.
Note that
%
\begin{eqnarray}\label{8}
|X\Delta|_2 &=& |Z\hat\theta-y-\Xi\hat\theta|_2\nonumber\\
&\le& \sqrt{n} (|Z\hat\theta-y|_\infty
+|\Xi\hat\theta|_\infty)\nonumber\\[-8pt]\\[-8pt]
&\le& \sqrt{n} (\delta|\hat\theta|_1
+|\Xi|_\infty|\hat\theta|_1 )\nonumber\\
&\le& 2\delta\sqrt{n} |\hat\theta|_1,\nonumber
\end{eqnarray}
which proves (\ref{6}).

Next, by the standard argument (cf. Lemma \ref{l1} below or, e.g.,
\cite{ct,brt}) we have
\[
|\Delta_{J^c}|_1\le|\Delta_{J}|_1.
\]
Thus,
%
\begin{equation}\label{7}
|\hat\theta-\theta_s|_1\le2|\Delta_{J}|_1 \le2 \sqrt{s}
|\Delta_{J}|_2 \le
\frac{2\sqrt{s}}{\kappa\sqrt{n}}|X\Delta|_2,
\end{equation}
where the last inequality follows from Assumption \ref{assumpREs}. Combining
(\ref{8}) and (\ref{7}), we get (\ref{4}).

To prove (\ref{5}), we introduce the set of indices $J_1$
corresponding to those $s$ coordinates of $\Delta$ outside
$J=J(\theta_s)$
which are largest in absolute value (we assume without loss of
generality that $2s\le p$). Define $J_{01}=J\cup J_1$. By a simple
argument that does not use any assumption (cf., e.g., \cite{ct,brt}
and the papers cited therein), we get
%
\begin{equation}\label{sim}
|\Delta_{J_{01}^c}|_2\le\frac{|\Delta_{J^c}|_1}{\sqrt{s}} .
\end{equation}
Thus,
\[
|\Delta_{J_{01}^c}|_2 \le\frac{|\Delta_{J}|_1}{\sqrt{s}} \le
|\Delta_{J}|_2\le|\Delta_{J_{01}}|_2,
\]
so that
\[
|\Delta|_2\le2|\Delta_{J_{01}}|_2.
\]
Now, by Assumption \hyperlink{assumpRE2s}{RE($2s$)} and (\ref{8}),
\[
|\Delta_{J_{01}}|_2\le\frac{1}{\kappa\sqrt{n}}|X\Delta|_2 \le
\frac{2\delta}{\kappa}|\hat\theta|_1
\]
and, hence, (\ref{5}) follows.

We finally prove (\ref{6a}). Note first that
%
\begin{eqnarray}\label{8a}
|\Psi(\hat\theta-\theta_s)|_\infty&\equiv& \frac1{n}|X^TX(\hat
\theta-\theta
_s)|_\infty\nonumber\\
&=& \frac1{n} \max_{1\le j\le p} \bigl|\mathbf{x}_{(j)}^TX(\hat\theta-\theta_s) \bigr|\\
&\le& \frac1{n}|X(\hat\theta-\theta_s)|_2 \max_{1\le j\le p}
\bigl|\mathbf{x}_{(j)}\bigr|_2 = \frac1{\sqrt{n}}|X(\hat\theta-\theta
_s)|_2\nonumber,
\end{eqnarray}
where $\mathbf{x}_{(j)}$ denotes the $j$th column of $X$ and the last
equality uses the fact that $|\mathbf{x}_{(j)}|_2=\sqrt{n}$ since
$|\mathbf{x}_{(j)}|_2^2/n$ are the diagonal elements of $X^T X/n$.
Therefore, by~(\ref{8}),
%
\begin{equation}\label{8ab}
|\Psi(\hat\theta-\theta_s)|_\infty\le2\delta|\hat\theta|_1.
\end{equation}
Now,
since the $j$th component of $\Psi(\hat\theta-\theta_s)$ is
\[
\bigl(\Psi(\hat\theta-\theta_s)\bigr)_j= (\hat\theta_j-\theta_{sj})+
\sum_{i=1,i\neq
j}^p\Psi_{ij}(\hat\theta_i-\theta_{si}),
\]
where $\theta_{si}$ is the $i$th component of $\theta_s$, we obtain
%
\begin{equation}\label{ll}
|\hat\theta-\theta_s|_\infty\le2\delta|\hat\theta|_1 + \rho
|\hat\theta-\theta_s|_1.
\end{equation}
Recall that Assumption \ref{assumpC} with $\rho<(3\alpha s)^{-1}$ implies Assumption
\ref{assumpREs} with $\kappa=\sqrt{1-1/\alpha}$ (cf. \cite{brt} or Lemma 2 in
\cite{l}). Thus, we can apply (\ref{4}) with this value of $\kappa$
to bound $|\hat\theta-\theta_s|_1$ in (\ref{ll}), which finally yields
(\ref{6a}).
\end{pf}
\begin{remark}\label{Remark1}
We can replace $|\hat\theta|_1$ by $|\theta_s|_1$
in all the inequalities of Theorem \ref{t1}.
\end{remark}
\begin{remark}\label{Remark2}
It is straightforward to deduce a bound
for $|\hat\theta-\theta_s|_q$ for any $1\le q\le2$ from the bounds
(\ref{4}) and (\ref{5}), as it is done, for example, in \cite{brt}.
\end{remark}

Under the assumptions of part (iv) of Theorem \ref{t1}, we get
%
\begin{equation}\label{6aa}
|\hat\theta-\theta_s|_\infty< C_*(\alpha)\delta|\hat\theta|_1,
\end{equation}
where $C_*(\alpha)=2 (1+\frac{2}{3\sqrt{\alpha(\alpha-1)}} )$
is a
constant. Based on this, we can define the thresholded estimator
$\tilde\theta=(\tilde\theta_1,\ldots,\tilde\theta_p)$, where
%
\begin{equation}\label{A}
\tilde\theta_j = \hat\theta_j I\{|\hat\theta_j|>\tau\},\qquad
j=1,\ldots,p,
\end{equation}
with the data-dependent\vspace*{1pt} threshold $\tau=C_*(\alpha)\delta|\hat
\theta|_1$ for
some $\alpha>1$. Here $I\{\cdot\}$ denotes the indicator function. It
is useful to note that since the MU-selector $\hat\theta$ is, in
general, not unique, the thresholded estimator $\tilde\theta$ is
also not necessarily unique.

Our next result shows that the thresholded estimator $\tilde\theta$
recovers the sparsity pattern and, moreover, it recovers the signs of
the coordinates of $s$-sparse solution $\theta_s$ (this property is
sometimes called the sign consistency; cf. \cite{mb,zy,zh,l}). We
define
\[
\operatorname{sign} \theta= \cases{
-1, &\quad if $\theta<0$,\cr
0, &\quad if $\theta=0$,\cr
1, &\quad if $\theta>0$.}
\]
\begin{theorem}\label{t2} Assume that $\theta_s\in\Theta$ is an $s$-sparse
solution of $y=X\theta$, and that $\Theta\subseteq\{\theta\in
\mathbb{R}^p\dvtx
|\theta|_1\le a\}$ for some $a>0$. Let (\ref{2}) and Assumption \ref{assumpC} hold
with $\rho< (3\alpha s)^{-1}$ for some $\alpha>1$. If
%
\begin{equation}\label{LB}
{\min_{j\in J(\theta_s)}} |\theta_{sj}|> C_*(\alpha)\delta{a},
\end{equation}
then
%
\begin{equation}\label{LBB}
\operatorname{sign} \tilde\theta_j = \operatorname{sign}
\theta_{sj},\qquad
j=1,\ldots,p,
\end{equation}
for all $\tilde\theta_j$ in (\ref{A}) such that $\hat\theta$ is an
MU-selector defined in (\ref{3}).
\end{theorem}
\begin{pf}
For $j\notin J(\theta_s)$ we have $\theta_{sj}=0$ and,
thus, by (\ref{6aa}), $|\hat\theta_j| = |\hat\theta_j-\theta
_{sj}|< C_*(\alpha
)\delta
|\hat\theta|_1= \tau$. Therefore, $\tilde\theta_j=0$ for $j\notin
J(\theta_s)$.
For $j\in J(\theta_s)$ note that (\ref{6aa}) implies
$|\hat\theta_j-\theta_{sj}|< C_*(\alpha)\delta{a}$. This and
assumption (\ref{LB})
yield that $\hat\theta_j$ has the same sign as $\theta_{sj}$.
\end{pf}
\begin{remark}\label{Remark3}
Note that, under Assumption \ref{assumpC} with $\rho<
(3\alpha s)^{-1}$ as required in Theorem \ref{t2}, the $s$-sparse
solution is unique; cf., for example, \cite{l}, page 93, so that the
right-hand side of (\ref{LBB}) is uniquely defined. The estimator
$\tilde\theta$ is not necessarily unique, nevertheless, Theorem \ref{t2}
assures that the sign recovery property (\ref{LBB}) holds for all
versions of $\tilde\theta$.
\end{remark}

\section{Sparse recovery for regression model with unknown design matrix}\label{sec4}

We consider now the general model (\ref{0}) and (\ref{00}) and assume
that it holds with an $s$-sparse vector of unknown parameters
$\theta^*=\theta_s\in\Theta$. Because of the presence of noise
$\xi$ that is
typically not small, we need to change the definition of
the MU-selector. We now define the MU-selector $\hat\theta$ as a
solution of the minimization problem
%
\begin{equation}\label{9}
\min\biggl\{ |\theta|_1\dvtx \theta\in\Theta,
\biggl|\frac1{n}Z^T(y-Z\theta) \biggr|_\infty\le(1+\delta)\delta|\theta|_1
+\varepsilon\biggr\}.
\end{equation}
Note that if $\delta=0$ and $\Theta=\mathbb{R}^p$, this
MU-selector becomes the
Dantzig selector of~\cite{ct}.

Similarly to (\ref{3}), the problem (\ref{9}) is a convex
minimization problem and it reduces to linear programming if
$\Theta=\mathbb{R}^p$, $\Theta$ is a linear subspace of $\mathbb
{R}^p$ or a simplex.

Throughout this section we will assume for simplicity that the
matrix $X$ is normalized, so that all the diagonal elements of the
Gram matrix $\Psi=X^TX/n$ are equal to 1. Extensions to general
matrices are straightforward, it only modifies the constants in
the expression $(1+\delta)\delta|\theta|_1 +\varepsilon$ in (\ref
{9}) and in the
theorems.

Note that under assumptions (\ref{1}) and (\ref{2}), the feasible set
of the convex problem (\ref{9}) is nonempty:
\[
\Theta_2 \equiv\biggl\{ \theta\in\Theta\dvtx
\biggl|\frac1{n}Z^T(y-Z\theta) \biggr|_\infty\le(1+\delta)\delta|\theta|_1
+\varepsilon\biggr\}
\ne
\varnothing.
\]
To prove this, let us show that the true vector $\theta^*=\theta_s$ belongs
to $\Theta_2$. In fact, by~(\ref{1}),
%
\begin{eqnarray}\label{9a}
\biggl|\frac1{n}Z^T(y-Z\theta_s) \biggr|_\infty&=&
\biggl|\frac1{n}Z^T(X\theta_s+\xi-Z\theta_s) \biggr|_\infty\nonumber\\
&\le& \biggl|\frac1{n}Z^T\xi\biggr|_\infty+
\biggl|\frac1{n}Z^T\Xi\theta_s \biggr|_\infty\\
&\le& \varepsilon+ \biggl|\frac1{n}Z^T\Xi\theta_s \biggr|_\infty.\nonumber
\end{eqnarray}
Next, note that, by (\ref{2}) and by the fact that all the diagonal elements
of $X^TX/n$ are equal to 1, the columns $\mathbf{z}_{(j)}$ of matrix
$Z$ satisfy $|\mathbf{z}_{(j)}|_2\le\sqrt{n}(1+\delta)$. Therefore,
arguing as in (\ref{8a}), we obtain
%
\begin{equation}\label{vaz}
\biggl|\frac1{n}Z^T\Xi\theta_s \biggr|_\infty\le
\frac{1+\delta}{\sqrt{n}} |\Xi\theta_s |_2 \le
(1+\delta) |\Xi\theta_s |_\infty\le(1+\delta)\delta|\theta_s|_1.
\end{equation}
This and (\ref{9a}) yield
\[
\biggl|\frac1{n}Z^T(y-Z\theta_s) \biggr|_\infty\le(1+\delta)\delta|\theta
_s|_1 +\varepsilon.
\]
Since we also assume that $\theta_s$ belongs to $\Theta$, the fact that
$\theta_s\in\Theta_2$ is proved. Thus, there always exists a solution
$\hat\theta$ of (\ref{9}). Of course, it is not necessarily unique.
%
%
\begin{theorem}\label{t3}
Assume that model (\ref{0}) and (\ref{00}) holds with an (unknown)
$s$-sparse parameter vector $\theta^*=\theta_s\in\Theta$ and
that all the
diagonal elements of $X^TX/n$ are equal to 1. Let (\ref{1}) and
(\ref{2}) hold. Set
\[
\nu= 2(2+\delta)\delta|\theta_s|_1 +2\varepsilon.
\]
Then for any solution $\hat\theta$ of (\ref{9}) we have the following
inequalities:

\begin{longlist}
\item
Under Assumption \ref{assumpREs}:
%
\begin{eqnarray}
\label{11a}
|\hat\theta-\theta_s|_1&\le& \frac{4\nu s}{\kappa^2} ,
\\
\label{10}
\frac{1}{n}|X(\hat\theta-\theta_s)|_2^2&\le& \frac{4\nu^2
s}{\kappa^2} .
\end{eqnarray}
\item Under Assumption \textup{RE($2s$)}:
%
\begin{equation}\label{11}
|\hat\theta-\theta_s|_q^q\le \biggl(\frac{4\nu}{\kappa^2} \biggr)^q
s\qquad
\forall1< q\le2.
\end{equation}
\item Under Assumption \ref{assumpC} with $\rho<\frac1{3\alpha s}, \alpha>1$:
%
\begin{equation}\label{14}
|\hat\theta-\theta_s|_\infty<
\frac{3\alpha+1}{3(\alpha-1)} \nu.
\end{equation}
\end{longlist}
\end{theorem}
\begin{pf}
Set $\Delta=\hat\theta-\theta_s$ and $J=J(\theta
_s)$. Note
first that (\ref{1}) and the fact that $\hat\theta$ belongs to the
feasible set $\Theta_2$ of (\ref{9}) imply
%
\begin{eqnarray}\label{xx}\qquad
\biggl|\frac1{n}X^TX\Delta\biggr|_\infty&\le&
\biggl|\frac1{n}Z^T(y-Z\hat\theta) \biggr|_\infty+
\biggl|\frac1{n}\Xi^TX(\hat\theta-\theta_s) \biggr|_\infty\nonumber\\
&&{} + \biggl|\frac1{n}Z^T\xi\biggr|_\infty+
\biggl|\frac1{n}Z^T\Xi\hat\theta\biggr|_\infty\\
&\le& (1+\delta)\delta|\hat\theta|_1 + 2\varepsilon+
\biggl|\frac1{n}Z^T\Xi\hat\theta\biggr|_\infty
+ \biggl|\frac1{n}\Xi^TX(\hat\theta-\theta_s) \biggr|_\infty.\nonumber
\end{eqnarray}
Now,
%
\begin{equation}\label{xix}
|\Xi^TX |_\infty= \max_{1\le j,k\le p} \bigl|\xi_{(j)}^T\mathbf{x}_{(k)}\bigr|
\le\max_{1\le j,k\le p} \bigl|\xi_{(j)}\bigr|_2 \bigl|\mathbf{x}_{(k)}\bigr|_2
\le\delta{n},
\end{equation}
where $\xi_{(j)}$ are the columns of $\Xi$ and we used that $|\mathbf{x}_{(k)}|_2=\sqrt{n}$ by assumption on $X^TX/n$, and
$|\xi_{(j)}|_2\le\delta\sqrt{n}$ by (\ref{2}). This implies
%
\begin{equation}\label{xixt}
\biggl|\frac1{n}\Xi^TX(\hat\theta-\theta_s) \biggr|_\infty\le|\hat\theta
-\theta_s|_1
\biggl|\frac1{n}\Xi^TX \biggr|_\infty\le\delta|\hat\theta-\theta_s|_1.
\end{equation}
Next, as in (\ref{vaz}), we obtain
%
\begin{equation}\label{xixt2}
\biggl|\frac1{n}Z^T\Xi\hat\theta\biggr|_\infty\le(1+\delta)\delta|\hat
\theta|_1.
\end{equation}
We now combine (\ref{xx}), (\ref{xixt}) and (\ref{xixt2}) to get
%
\begin{equation}\label{xixt3}
\biggl|\frac1{n}X^TX\Delta\biggr|_\infty\le2\varepsilon+ 2(1+\delta)\delta
|\hat\theta|_1+
\delta|\hat\theta-\theta_s|_1 \le\nu.
\end{equation}
Taking into account (\ref{xixt3}), the proof of (\ref{11a}),
(\ref{10}) and (\ref{11}) follows the same lines as the proof of
Theorem 7.1 in \cite{brt} where we should set $r=\nu/2$, $m=s$.

We now prove (\ref{14}). We proceed as in the proof of (\ref{6a}) in
Theorem \ref{t1}, with the only difference that now we replace
(\ref{8ab}) by (\ref{xixt3}). Thus, instead of (\ref{ll}), we obtain
%
\begin{equation}\label{lo}
|\hat\theta-\theta_s|_\infty\le\nu+ \rho|\hat\theta-\theta_s|_1.
\end{equation}
Next, recall that Assumption \ref{assumpC} with $\rho<(3\alpha s)^{-1}$ implies
Assumption \ref{assumpREs} with $\kappa^2 = 1-1/\alpha$ (cf. \cite{brt} or Lemma
2 in \cite{l}). Using in (\ref{lo}) the bound (\ref{11}) with $q=1$
and $\kappa^2 = 1-1/\alpha$, we obtain (\ref{14}). This finishes the
proof of the theorem.
\end{pf}

Note that, in contrast to Theorem \ref{t1}, the bounds of Theorem
\ref{t3} do not depend on $|\hat\theta|_1$ but on the unknown
$|\theta_s|_1$ (cf. definition of $\nu$). This drawback can be
corrected for small values of $\delta$, as shown in the next result.
%
%
\begin{theorem}\label{t3_1}
Let the assumptions of Theorem \ref{t3} hold and
$\delta<\frac{\kappa^2}{4s}$. Set
\[
\nu_1 = 2(1+\delta)\delta|\hat\theta|_1 +2\varepsilon.
\]
Then for any solution $\hat\theta$ of (\ref{9}) we have the following
inequalities:

\begin{longlist}
\item
Under Assumption \ref{assumpREs}:
%
\begin{eqnarray}\label{11ab}
|\hat\theta-\theta_s|_1&\le& \frac{4\nu_1 s}{\kappa^2}
\biggl(1-\frac
{4\delta{s}}{\kappa^2} \biggr)^{-1} ,
\\
\label{10b}
\frac{1}{n}|X(\hat\theta-\theta_s)|_2^2&\le& \frac{4\nu_1^2
s}{\kappa^2}\biggl(1-\frac{4\delta{s}}{\kappa^2} \biggr)^{-2} .
\end{eqnarray}
\item Under Assumption \textup{RE($2s$)}:
%
\begin{equation}\label{11b}
|\hat\theta-\theta_s|_q^q\le
\biggl(\frac{4\nu_1}{\kappa^2} \biggr)^q \biggl(1-\frac{4\delta{s}}{\kappa^2}
\biggr)^{-q} s\qquad \forall1< q\le2.
\end{equation}
\item Under Assumption \ref{assumpC} with $\rho<\frac1{3\alpha s}, \alpha>1$, and
$\delta\le\frac{\kappa^2}{8s}$:
%
\begin{equation}\label{14b}
|\hat\theta-\theta_s|_\infty<
\frac{2(3\alpha+1)}{3(\alpha-1)} \nu_1.
\end{equation}
\end{longlist}
\end{theorem}
\begin{pf}
We use the same notation as in the proof of
Theorem \ref{t3}. From (\ref{xixt3}) and the fact that
$|\Delta_{J^c}|_1\le|\Delta_{J}|_1$, we obtain
%
\begin{eqnarray}\label{t3_11}
\frac{1}{n}|X\Delta|_2^2&\le&
|\Delta|_1 \biggl|\frac1{n}X^TX\Delta\biggr|_\infty\nonumber\\
&\le& |\Delta|_1(\nu_1+ \delta|\Delta|_1 ) \nonumber\\[-8pt]\\[-8pt]
&\le& 2|\Delta_J|_1(\nu_1+ 2\delta|\Delta_J|_1 ) \nonumber
\\
&\le& 2\sqrt{s}|\Delta_J|_2\bigl(\nu_1+ 2\sqrt{s}\delta|\Delta_J|_2 \bigr).
\nonumber
\end{eqnarray}
Similar arguments as for (\ref{7}) easily yield the inequality
%
\begin{equation}\label{t3_12}
|\Delta_J|_2 \le\frac{2\sqrt{s}\nu_1}{\kappa^2} \biggl(1-\frac
{4\delta
{s}}{\kappa^2} \biggr)^{-1}
\end{equation}
and (\ref{11ab}). In the same way, (\ref{11b}) deduces from
(\ref{t3_12}) following the analogous part of the proof of Theorem
7.1 in \cite{brt} where we should set $r=(\nu_1/2) (1-4\delta
{s}/\kappa^2 )^{-1}$ and $m=s$.

Finally, to get the sup-norm inequality (\ref{14b}), we proceed as in
the proof of (\ref{6a}) in Theorem \ref{t1} or in that of (\ref{14})
in Theorem \ref{t3}, with the only difference that instead of
(\ref{ll}) we use the bound
\[
|\hat\theta-\theta_s|_\infty\le\nu_1 + (\rho+\delta) |\hat
\theta-\theta_s|_1
\]
that follows from (\ref{xixt3}) and the fact that
$\delta\le\frac{\kappa^2}{8s}$. This finishes the proof of
Theorem~\ref{t3_1}.
\end{pf}

As in Section \ref{sec3}, we now define a
thresholded estimator $\tilde\theta=(\tilde\theta_1,\ldots,\tilde
\theta_p)$ by
the formula
%
\begin{equation}\label{B}
\tilde\theta_j = \hat\theta_j I\{|\hat\theta_j|>\tau_1\},\qquad
j=1,\ldots,p,
\end{equation}
where the threshold is given either by
%
\begin{equation}\label{B_1}
\tau_1=\frac{3\alpha+1}{3(\alpha-1)}
\bigl(2\varepsilon+
2(2+\delta)\delta{a} \bigr)
\end{equation}
for $\alpha>1$, $a>0$, or by
%
\begin{equation}\label{B_2}
\tau_1=\frac{2(3\alpha+1)}{3(\alpha
-1)} \bigl(2\varepsilon+
2(1+\delta)\delta|\hat\theta|_1 \bigr)
\end{equation}
for $\alpha>1$. Note that the threshold (\ref{B_2}) is completely
data-driven if $\varepsilon$ and $\delta$ are known.

The next theorem shows that under some assumptions the thresholded
estimator defined in (\ref{B}) recovers the sparsity pattern and,
moreover, it recovers the signs of the coordinates of the $s$-sparse
solution $\theta_s$.
\begin{theorem}\label{t4}
Assume that model (\ref{0}) and (\ref{00}) holds with the $s$-sparse vector
of unknown parameters\vspace*{1pt} $\theta^*=\theta_s\in\Theta$ and that
(\ref{1}),
(\ref{2}) and Assumption \ref{assumpC} hold with $\rho< (3\alpha s)^{-1}$ for some
$\alpha>1$. Let either $\Theta\subseteq\{\theta\in\mathbb
{R}^p\dvtx |\theta|_1\le a\}$ for
some $a>0$ and the threshold $\tau_1$ is given by (\ref{B_1}), or
$\delta\le\frac{\kappa^2}{8s}$ and the threshold $\tau_1$ is
given by~(\ref{B_2}). If
%
\begin{equation}\label{LB1}
{\min_{j\in J(\theta_s)}} |\theta_{sj}|> \tau_1,
\end{equation}
then
%
\begin{equation}\label{LBB1}
\operatorname{sign} \tilde\theta_j = \operatorname{sign}
\theta_{sj},\qquad
j=1,\ldots,p,
\end{equation}
for all $\tilde\theta_j$ in (\ref{B}) where $\hat\theta$ is a
MU-selector defined in (\ref{9}).
\end{theorem}
\begin{pf}
It goes along the same lines as the proof of
Theorem \ref{t2}.

We can make here the same remarks as in Section \ref{sec3} about the
nonuniqueness of the estimators. Indeed, $\tilde\theta$ is not
necessarily unique, but Theorem \ref{t4} assures the sign
recovery property (\ref{LBB1}) holds for all versions of $\tilde
\theta$.
\end{pf}
\begin{remark}\label{Remark4}
The argument of this section can be applied with
minor modifications to the model
\begin{eqnarray*}
y&=& Z\theta^* + \xi,\\
Z&=&X+\Xi.
\end{eqnarray*}
This is no longer the errors-in-variables setting, but just the
usual regression setting where $X$ is some ``nominal'' design matrix
and $\Xi$ can be viewed as its perturbation. The results of this
section suggest that small perturbations of the design matrix $X$
beyond the restricted eigenvalue condition are in fact not so
dangerous, once the method of recovery is chosen in a proper way.
Indeed, such perturbations lead to the extra terms in the bounds
proportional to the $\ell_1$-norm of the solution. Roughly speaking,
our bounds suggest that the MU-selector is robust with respect to
possible violations of the the restricted eigenvalue condition,
provided that the perturbations are small enough and the
$\ell_1$-norm of the true $\theta$ is reasonably bounded. This
offers a possible way of relaxing the strong conditions usually
imposed in the context of $\ell_1$-penalized sparse estimation. Note
that another way to do it can be found in \cite{dt1,dt2}, suggesting
a computationally feasible method of sparse estimation with no
assumption on $X$. However, the oracle inequalities of
\cite{dt1,dt2} hold only for the prediction risk.
\end{remark}

\section{Approximately $s$-sparse solutions}\label{sec5}

The results of the previous sections can be easily generalized to
the setting where the true $\theta^*$ is arbitrary, not necessarily
$s$-sparse. This might be of interest in the context of inverse
problems with unknown operator, as discussed in the \hyperref[intro]{Introduction}.
Then the bounds will involve a residual term, which is a
difference between $\theta^*$ and its $s$-sparse approximation~$\theta_s$.
In particular, we can take $\theta_s$ as the best $s$-sparse
approximation of $\theta^*$, that is, the vector that coincides with
$\theta^*$ in the $s$ coordinates with largest absolute values and has
other coordinates that vanish.

We will use the following slightly strengthened version of
Assumption \ref{assumpREs}, where we only increase a numerical constant in
the definition of the set over which the minimum is taken (cf.
\cite{brt}).
\renewcommand{\theassumption}{RE($s,2$)}
\begin{assumption}\label{assumpREs2}
There exists $\kappa>0$ such that
\[
\min_{\Delta\ne0\dvtx |\Delta_{J^c}|_1\le
2|\Delta_{J}|_1}\frac{|X\Delta|_2}{\sqrt{n}|\Delta_J|_2} \ge
\kappa
\]
for all subsets $J$ of $\{1,\ldots,p\}$ of cardinality $|J|\le s$.
\end{assumption}

It is easy to check that Assumption \ref{assumpC} with $\rho<\frac1{5\alpha s}$ for
some $\alpha>1$ implies Assumption \ref{assumpREs2} with $\kappa^2=
1-1/\alpha$
(cf. \cite{brt}).

We now state the main result of this section.
\begin{theorem}\label{t5}
Assume that there exists a solution $\theta^*\in\Theta$ of the equation
$y=X\theta$. Let (\ref{2}) hold. Then for any solution $\hat\theta$
of (\ref{3}) we have the following inequalities:

\begin{longlist}
\item
%
\begin{equation}\label{6x}
\frac{1}{n}|X(\hat\theta-\theta^*)|_2^2\le 4\delta^2|\hat
\theta
|_1^2.
\end{equation}
\item If Assumption \ref{assumpREs2} holds, then
%
\begin{equation}\label{4x}
|\hat\theta-\theta^*|_1\le \frac{4\sqrt{s}\delta}{\kappa
}|\hat\theta
|_1 +
6\min_{J\dvtx |J|\le s} |\theta^*_{J^c}|_1.
\end{equation}
\item If Assumption \ref{assumpC} holds with $\rho<\frac1{5\alpha s}, \alpha
>1$, then
%
\begin{equation}\label{6ax}
|\hat\theta-\theta^*|_\infty<
2 \biggl(1+\frac{2}{5\sqrt{s\alpha(\alpha-1)}} \biggr)\delta|\hat\theta
|_1 +
\frac{6}{5\alpha s} \min_{J\dvtx |J|\le s} |\theta^*_{J^c}|_1.
\end{equation}
\end{longlist}
\end{theorem}
\begin{pf}
Set $\Delta=\hat\theta-\theta^*$ and let
$J\subset\{1,\ldots,p\}$ be an arbitrary set of indices such that
$|J|\le s$. First, note that (\ref{6x}) is already proved in
Theorem \ref{t1}, since (\ref{8}) is valid with $\Delta=\hat
\theta-\theta^*$.

We will use the following
elementary fact (cf., e.g., \cite{c08,crt,ct1,ct}) that we state for
convenience as a lemma.
\begin{lemma}\label{l1}
Let $\hat\theta$ be a solution of the problem
\[
\min\{ |\theta|_1\dvtx \theta\in\Theta'\},
\]
where $\Theta'$ is a subset of $\mathbb{R}^p$. Let $\theta^*$ be
any element of
$\Theta'$ and $J$ any subset of $\{1,\ldots,p\}$. Then for
$\Delta=\hat\theta-\theta^*$ we have
%
\begin{equation}\label{l11}
|\Delta_{J^c}|_1\le|\Delta_{J}|_1 +
2|\theta^*_{J^c}|_1.
\end{equation}
\end{lemma}
\begin{pf}
\begin{eqnarray*}
|\theta^*_{J}|_1+|\theta^*_{J^c}|_1 &=& |\theta^*|_1 \ge|\hat
\theta|_1
=|\hat\theta_J|_1
+ |\hat\theta_{J^c}|_1\\
&=&|\Delta_J + \theta^*_{J} |_1 + |\Delta_{J^c}+\theta
^*_{J^c}|_1 \\&\ge& |\theta^*_{J}
|_1 - |\Delta_J|_1 + |\Delta_{J^c}|_1 - |\theta^*_{J^c}|_1.
\end{eqnarray*}
\upqed\end{pf}

To prove (\ref{4x}), consider separately the following two cases:
(a) $2|\theta^*_{J^c}|_1\le|\Delta_J|_1$ and (b) $2|\theta^*_{J^c}|_1>
|\Delta_J|_1$. In case (a) we use (\ref{l11}) to obtain
$|\Delta_{J^c}|_1\le2|\Delta_J|_1$. Therefore, by Assumption \ref{assumpREs2}
and (\ref{6x}),
\[
|\Delta_{J}|_2 \le\frac1{\kappa\sqrt{n}}|X\Delta|_2
\le\frac{2\delta}{\kappa}|\hat\theta|_1.
\]
This and (\ref{l11}) imply that, in case (a),
%
\begin{eqnarray}\label{5_1}
|\Delta|_1&\le& 2|\Delta_{J}|_1 + 2|\theta^*_{J^c}|_1 \le2\sqrt
{s}|\Delta_{J}|_2+
2|\theta^*_{J^c}|_1 \nonumber\\[-8pt]\\[-8pt]
&\le&
\frac{4\sqrt{s}\delta}{\kappa}|\hat\theta|_1+2|\theta^*_{J^c}|_1.\nonumber
\end{eqnarray}
In case (b) we immediately deduce from (\ref{l11}) that $|\Delta
|_1\le
6 |\theta^*_{J^c}|_1$. Combining this with (\ref{5_1}), we obtain
(\ref{4x}).

To prove (\ref{6ax}), note that the argument leading to (\ref{ll}) is
applicable here with $\theta^*$ in place of $\theta_s$. Thus,
%
\begin{equation}\label{lll}
|\hat\theta-\theta^*|_\infty\le2\delta|\hat\theta|_1 + \rho
|\hat\theta-\theta^*|_1.
\end{equation}
Now, as mentioned above, Assumption \ref{assumpC} with $\rho<\frac1{5\alpha s}$,
$\alpha>1$, implies Assumption \ref{assumpREs2} with $\kappa^2= 1-1/\alpha
$. Using
(\ref{4x}) with this value of $\kappa$ to bound $|\hat\theta
-\theta_s|_1$ in
(\ref{lll}), we arrive at (\ref{6ax}). This proves the theorem.
\end{pf}

Note that under Assumption \ref{assumpC} we can also bound the $\ell_2$ norm
of the difference $\hat\theta-\theta^*$, as well as all its $\ell_r$
norms with $r>1$. However, Assumption \ref{assumpC} is rather restrictive.
For instance, it is not valid for Toeplitz matrices $\Psi$ or for
matrices $X$ with independent standard Gaussian entries (for the
latter case, Assumption RE is assured with overwhelming
probability if $s$ is of a smaller order than $n/\log p$). The
next theorem shows that we can bound correctly the $\ell_2$ norm
$|\hat\theta-\theta^*|_2$ under the following condition which is weaker
than Assumption \ref{assumpC} but somewhat stronger than Assumption RE.
\renewcommand{\theassumption}{RE$^{\prime}$($s,2$)}
\begin{assumption}\label{assumpREprimes2}\hypertarget{aaa}
There exist $\kappa>0$
and $c_1>0$ such that
%
\begin{equation}
\label{reprime}\min_{\Delta\ne0\dvtx |\Delta_{J^c}|_1\le
2|\Delta_{J}|_1+a}\frac{|X\Delta|_2^2/n + c_1 a
|\Delta_J|_2/\sqrt{s}}{|\Delta_J|_2^2} \ge\kappa^2
\end{equation}
for all $a\ge0$ and all subsets $J$ of $\{1,\ldots,p\}$ of
cardinality $|J|\le s$.
\end{assumption}

Note that Assumption \ref{assumpREs2} is a special case of
(\ref{reprime}) corresponding to $a=0$. Note also that Assumption
\ref{assumpREprimes2} is satisfied if the restricted isometry
assumption \cite{crt,ct1,ct} holds with the isometry coefficient
close enough to 1. This is not hard to show following the lines of
\cite{c08}.
\begin{theorem}\label{t6}
Assume that\vspace*{1pt} there exists a solution $\theta^*\in\Theta$ of the equation
$y=X\theta$. Let (\ref{2}) and Assumption \textup{\hyperlink{aaa}{RE$^{\prime}$($2s,2$)}} hold.
Then for any solution $\hat\theta$ of (\ref{3}) we have
%
\begin{equation}\label{5x}
|\hat\theta-\theta^*|_2\le \frac{4\delta}{\kappa}|\hat
\theta|_1 + \biggl(4 +
\frac{2\sqrt{c_1}}{\kappa} \biggr) \min_{J\dvtx |J|\le s}
\frac{|\theta^*_{J^c}|_1}{\sqrt{s}}.
\end{equation}
\end{theorem}
\begin{pf}
Set, as before, $\Delta=\hat\theta-\theta^*$ and let
$J\subset\{1,\ldots,p\}$ be an arbitrary set of indices such that
$|J|\le s$.
We first note that (\ref{l11}), (\ref{sim}) and the fact
that $|\Delta_{J}|_1\le\sqrt{s}|\Delta_{J}|_2$ imply
%
\begin{equation}\label{5_2}
|\Delta_{J_{01}^c}|_2\le|\Delta_{J}|_2 + \frac{2}{\sqrt
{s}}|\theta^*_{J^c}|_1.
\end{equation}
Consider separately the cases $2|\theta^*_{J^c}|_1/\sqrt{s}\le
|\Delta_{J_{01}}|_2$ and $2|\theta^*_{J^c}|_1/\sqrt{s}> |\Delta
_{J_{01}}|_2$.

(a) In the case $2|\theta^*_{J^c}|_1/\sqrt{s}\le|\Delta
_{J_{01}}|_2$ we
have $|\Delta_{J_{01}^c}|_1\le2|\Delta_{J_{01}}|_1$. Also,
$|J_{01}|\le
2s$ by the definition of $J_{01}$. Therefore, using Assumption
\hyperlink{aaa}{RE$^{\prime}$($2s,2$)} with $a=0$ and~(\ref{6x}), we get
\[
|\Delta_{J_{01}}|_2 \le\frac1{\kappa\sqrt{n}}|X\Delta|_2
\le\frac{2\delta}{\kappa}|\hat\theta|_1.
\]
This and (\ref{5_2}) imply
%
\begin{equation}\label{5_3}
|\Delta|_2\le|\Delta_{J_{01}}|_2 + |\Delta_{J}|_2 +
\frac{2}{\sqrt{s}}|\theta^*_{J^c}|_1
\le\frac{4\delta}{\kappa}|\hat\theta|_1 + \frac{2}{\sqrt
{s}}|\theta
^*_{J^c}|_1.
\end{equation}
Thus, (\ref{5x}) is proved in the case $2|\theta^*_{J^c}|_1/\sqrt
{s}\le
|\Delta_{J_{01}}|_2$.

(b) It remains to prove (\ref{5x}) in the case
$2|\theta^*_{J^c}|_1/\sqrt{s}> |\Delta_{J_{01}}|_2$. This
condition and
(\ref{5_2}) immediately yield
\[
|\Delta_{J_{01}^c}|_2\le
\frac{4}{\sqrt{s}}|\theta^*_{J^c}|_1,
\]
so that
%
\begin{equation}\label{5_5}|\Delta|_2\le|\Delta_{J_{01}}|_2 +
\frac{4}{\sqrt{s}}|\theta^*_{J^c}|_1.
\end{equation}
Next, from (\ref{l11}) we easily get
\[
|\Delta_{J_{01}^c}|_1\le|\Delta
_{J_{01}}|_1 +
2|\theta^*_{J^c}|_1.
\]
Therefore, using Assumption \hyperlink{aaa}{RE$^{\prime}$($2s,2$)} with
$a=2|\theta^*_{J^c}|_1$ and (\ref{6x}), we find
\begin{eqnarray*}
\kappa^2|\Delta_{J_{01}}|_2^2&\le& \frac1{n}|X\Delta|_2^2 + 2 c_1
\frac{|\theta^*_{J^c}|_1|\Delta_{J_{01}}|_2} {\sqrt{2s}}\\
&\le& 4\delta^2|\hat\theta|_1^2+ 2\sqrt{2}c_1 \frac{|\theta
^*_{J^c}|_1^2}{s},
\end{eqnarray*}
where we used that $2|\theta^*_{J^c}|_1/\sqrt{s}> |\Delta
_{J_{01}}|_2$. The
last display and (\ref{5_5}) imply that (\ref{5x}) holds in the case
$2|\theta^*_{J^c}|_1/\sqrt{s}> |\Delta_{J_{01}}|_2$.
\end{pf}


\section{Random noise}\label{sec6}

If $\xi$ and $\Xi$ are random and
conditions (\ref{1}) and (\ref{2}) are satisfied with a probability
close to 1, then all the bounds in the above theorems remain valid
with the same probability. This holds in different situations under
natural assumptions that we briefly discuss in this section.

First, it is not hard to see that if $\xi$ is normal with zero mean
and covariance matrix $\sigma^2 I$ where $I$ denotes the identity
matrix, and we take
%
\begin{equation}\label{choice_e}
\varepsilon=A\sigma\sqrt{\frac{\log p}{n}}
\end{equation}
for some $A>(1+\delta)\sqrt{2}$, then condition (\ref{1}) holds with
probability at least $1-p^{1-A^2/2}$. If $p$ is very large, this
probability is very close to 1. A similar remark holds for sub-Gaussian
$\xi$.

For more general $\xi$ we can guarantee condition (\ref{1}) only
with a larger value of $\varepsilon$ and with a probability that is
not as
close to 1 as in the Gaussian case. For example, if the components
$\xi_i$ of $\xi$ are independent zero mean random variables with
uniformly bounded variances, $E(\xi_i^2)\le\sigma^2<\infty$,
$i=1,\ldots,n$, and if the elements $X_{ij}, i=1,\ldots,n, j=1,\ldots,
p,$ of matrix $X$ satisfy
\[
\frac1{n}\sum_{i=1}^n \max_{j=1,\ldots,p}|X_{ij}|^2 \le c
\]
for some constant $c$, then condition (\ref{1}) holds with
probability at least $1-O (\frac{\log p}{\varepsilon^2n} )$~\cite{l}.
In particular, we can take
\[
\varepsilon= A\sqrt{\frac{(\log p)^{1+\gamma}}{n}} ,
\]
and then condition (\ref{1}) holds with probability at least
$1-O ((\log p)^{-\gamma} )$.

For the choice of $\delta$ in condition (\ref{2}) we can consider
the examples related to portfolio selection and to inverse problems
with unknown operator; cf. Section~\ref{sec2}. In both examples we have
repeated measurements. The matrix $Z$ is either the average of
several observed matrices with mean $X$, or the empirical covariance
matrix, with $X$ defined as the corresponding population covariance
matrix (in the latter case $p=n$). Then the threshold $\delta$ in
condition (\ref{2}) can be determined in the same spirit as
$\varepsilon$ in
condition (\ref{1}). We omit further details.

Finally, consider the model with missing data discussed in Section
\ref{sec2}. In this example direct application of condition (\ref{2}) leads
to bounds which are too loose. Indeed, $\delta$ can be of the order of
$|X|_\infty$. However, we argue that the MU-selector of the form
(\ref{mu_sel}) with suitable $\lambda$ still satisfies good bounds
if the probability $\pi$ that an entry of $X$ is not observed
remains small. This needs a refinement of our argument for the
particular setting. We sketch it now. Note first that under the
assumptions of Theorem \ref{t3} for a deterministic matrix $X$ and
for $Z_{ij}= X_{ij} + \xi_{ij}^\prime$, where $\xi_{ij}^\prime$
are defined in Section \ref{sec2}, we have, with probability close to 1 when
$n$ is large,
%
\begin{eqnarray}\label{new_ass}
\biggl|\frac1{n}\Xi^TX \biggr|_\infty&\le&\delta_1 ,\qquad
\biggl|\frac1{n}X^T\Xi\biggr|_\infty\le\delta_1,\\
\label{new_ass1}
\biggl|\frac1{n}\bigl(\Xi^T\Xi- \operatorname{diag}(\Xi^T\Xi)\bigr) \biggr|_\infty
&\le&\delta_2 , \\
\label{new_ass2}
\biggl|\frac1{n}\operatorname{diag}(\Xi^T\Xi) \biggr|_\infty&\le& C \pi,
\end{eqnarray}
where $\Xi$ is the matrix with entries $\xi_{ij}^\prime$,
$\operatorname{diag}(\Xi^T\Xi)$ denotes the diagonal matrix having the same
diagonal elements as $\Xi^T\Xi$, $C>0$ is a constant, and $\delta_1,
\delta_2>0$ are small if $n$ is large. Indeed, (\ref{new_ass}) and
(\ref{new_ass1}) follow from the standard properties of zero mean
sub-Gaussian variables, while (\ref{new_ass2}) is due to the fact
that the expectations of the diagonal elements of
$\frac1{n}\Xi^T\Xi$ are proportional to $\pi$.

We now observe that under assumptions
(\ref{new_ass}) and (\ref{new_ass1}) the constant $(1+\delta)\delta$ in (\ref{vaz})
can be
replaced by $\delta_1+\delta_2+C\pi$. This motivates the use of the
MU-selector (\ref{mu_sel}) with $\lambda=\delta_1+\delta_2+C\pi
$. For such an
MU-selector we have an analog of Theorem~\ref{t3} if we replace
assumption (\ref{2}) by assumptions (\ref{new_ass}) and (\ref{new_ass1}).
The only difference is in the form of $\nu$ which
now becomes a linear combination of $\delta_1, \delta_2$ and $\pi
$. This new
value of $\nu$ is small for $n$ large enough and small $\pi$. In
conclusion, the MU-selector (\ref{mu_sel}) with suitable $\lambda$
achieves good theoretical bounds provided that $\pi$ is small enough
and $n$ is large. This is confirmed by simulations in the next
section.


\section{Numerical experiments}\label{appli}
We present here three illustrative numerical applications. The
first two are based on simulated data and the last one on real
data.
\subsection{Censored matrix}\label{cens}
We begin with a model where we only observe censored elements of
the matrix $X$. More precisely, for a positive censoring value
$t$, instead of $X_{ij}$, we observe
%
\begin{equation}\label{censored}
Z_{ij}=X_{ij}I\{|X_{ij}|\leq
t\}+t(\operatorname{sign}X_{ij})I\{|X_{ij}|>t\}.
\end{equation}

\subsubsection*{Experiment}

\mbox{}

-- We take a matrix $X$ of size $100 \times500$
($n=100, p=500$) which is the normalized version (centered and
then normalized so that all the diagonal elements of the
associated Gram matrix are equal to 1) of a $100
\times500$ matrix with i.i.d. standard Gaussian entries.

-- For a given integer $s$, we randomly (uniformly)
choose $s$ nonzero elements in a vector $\theta$ of size $500$.
The associated values are equal to $0.5$. We will take $s=1,2,3,5,10$.

-- We set $y=X\theta+\xi$, where $\xi$ is a normal random
vector with zero mean
and covariance matrix $\sigma^2 I$ where $\sigma=0.05/1.96$ (so that
for an element of $\xi$, the probability of being between $-0.05$ and
$0.05$ is $95\%$).

-- We compute the matrix $Z$ following (\ref{censored})
with
$t=0.9.$

-- We run a linear programming algorithm to compute the
solution of (\ref{9}) where we optimize over $\Theta
=\mathbb{R}_{+}^{500}$. The value of $\varepsilon$ is chosen following
(\ref{choice_e}) with $A=(1+\delta)\sqrt{2}$. We note here that in the
simulations below the choice of $\varepsilon$ is not crucial because the
terms with $\delta$ in the definition of the estimator are of a larger
order of magnitude. Varying $\varepsilon$ within a sufficiently wide range
does not essentially modify the simulation results. The choice of
parameter $\delta$ is done the following way.

\subsubsection*{Choice of $\delta$}
\begin{figure}

\includegraphics{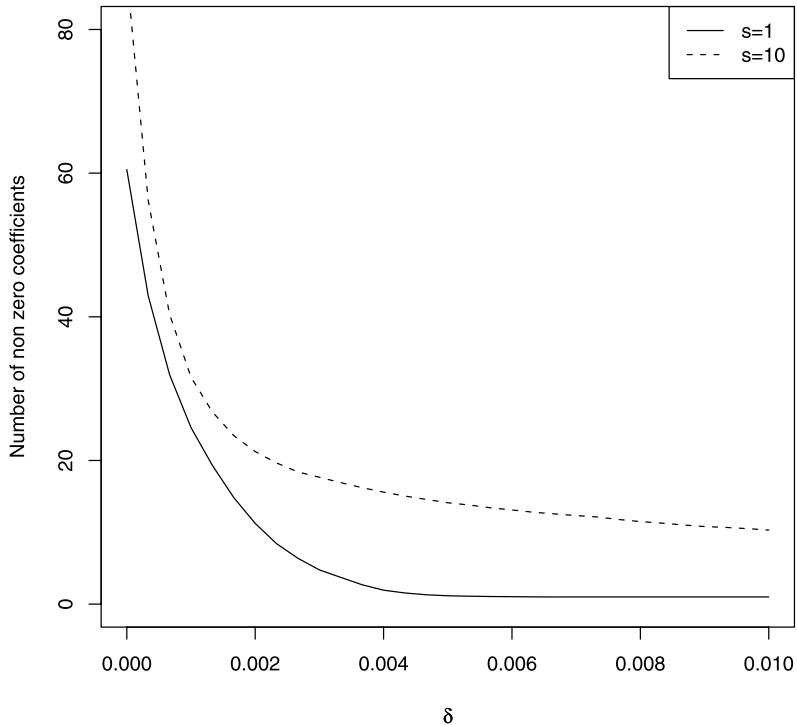}

\caption{Average number of nonzero coefficients in the model with
censored matrix for $s=1$ and $s=10$.}
\label{figure1}
\end{figure}

The choice of $\delta$ in practice is quite crucial. A very small
value of $\delta$ means that the matrix uncertainty is not taken
into account, whereas a too large value of $\delta$ means that we
overestimate this uncertainty. In both situations the resulting
estimator exhibits poor behavior. Consequently, in practice, it is
important to select $\delta$ within a reasonable range of values.
We suggest to choose the range of candidate $\delta$ with the
``elbow'' rule. We plot the number of retrieved nonzero
coefficients as a function of $\delta$. Then we consider that a
value of $\delta$ can be chosen only if the plot is (or begins to
be) flat around it. Usually such a plot is highly decreasing at the
beginning and then stabilizes; cf. Figure \ref{figure1}. Following this, we take
the values in the flat zone $\delta=0.05,0.75,0.1$ for $s=1,2,3,5$
and $\delta=0.01,0.05,0.1$ for $s=10$ (the plot for $s=10$ suggests
to start with smaller values for $\delta$).

-- We also compute the Lasso estimator with Mallows'
$C_p$ choice of the tuning parameter (we use the Lars R-package of
T. Hastie and B. Efron) and the Dantzig selector of \cite{ct},
with the same value $\varepsilon$. Moreover, we compute the thresholded
versions of the estimators (T-Lasso, T-Dantzig, T-$\delta$). More
precisely, the retrieved coefficients whose absolute values are
smaller than $20\%$ of the true value of the nonzero coefficients
(i.e., smaller than 0.1) are set to zero.

-- For all the considered estimators $\hat{\theta}$ of
$\theta$ we compute the error measures
\[
\operatorname{Err}_1=|\hat{\theta}-\theta|_2^2\quad\mbox{and}\quad
\operatorname{Err}_2=
|X(\hat{\theta}-\theta)|_2^2.
\]
We also record the retrieved
sparsity pattern, which is defined as the set of the nonzero
coefficients of $\hat{\theta}$.

-- For each value of $s$ we run $100$ Monte Carlo
simulations.

\subsubsection*{Results}
Tables \ref{table1}--\ref{table5} present the empirical averages and standard deviations
(in brackets) of $\operatorname{Err}_1$, $\operatorname{Err}_2$, of
the number of
nonzero coefficients in $\hat{\theta}$ ($\operatorname{Nb}_1$) and of the number
of nonzero coefficients in $\hat{\theta}$ belonging to the true
sparsity pattern ($\operatorname{Nb}_2$). We also present the total number of
simulations where the sparsity pattern is exactly retrieved (Exact).
Note that here and in the next numerical examples when a coefficient
belonging to the sparsity pattern is retrieved it has
systematically the correct sign.

\begin{table}
\caption{Results for the model with censored matrix, $s=1$}\label{table1}
{\fontsize{8.8}{10.65}\selectfont{\begin{tabular*}{\tablewidth}{@{\extracolsep{\fill}}lcd{2.5}ccd{3.0}@{}}
\hline
& \multicolumn{1}{c}{$\bolds{\operatorname{Err}_1}$}
& \multicolumn{1}{c}{$\bolds{\operatorname{Err}_2}$}
& \multicolumn{1}{c}{$\bolds{\operatorname{Nb}_1}$}
& \multicolumn{1}{c}{$\bolds{\operatorname{Nb}_2}$}
& \multicolumn{1}{c@{}}{\textbf{Exact}}\\
\hline
Lasso &0.0679&12.33&\multicolumn{1}{r}{95.20\phantom{0)}}&1&0\\
&(0.0128)&(2.016)&\multicolumn{1}{r}{(2.245)}&(0)&\\
T-Lasso&0.0271&2.712&1&1&100\\
&(0.0086)&(0.8615)&(0)&(0)&\\
Dantzig &0.0399&3.982&\multicolumn{1}{r}{56.92\phantom{0)}}&1&0\\
&(0.0076)&(0.9880)&\multicolumn{1}{r}{(5.594)}&(0)&\\
T-Dantzig&0.0260&2.599&1&1&100\\
&(0.0068)&(0.6860)&(0)&(0)&\\
$\delta=0.05$ &0.0122&1.231&\multicolumn{1}{r}{1.16\phantom{0)}}&1&85\\
&(0.0027)&(0.2783)&\multicolumn{1}{r}{(0.393)}&(0)&\\
T-$\delta=0.05$&0.0122&1.224&1&1&100\\
&(0.0028)&(0.2816)&(0)&(0)&\\
$\delta=0.075$ &0.0064&0.649&1&1&100\\
&(0.0017)&(0.1715)&(0)&(0)&\\
T-$\delta=0.075$&0.0064&0.649&1&1&100\\
&(0.0017)&(0.1715)&(0)&(0)&\\
$\delta=0.1$&0.0023&0.2330&1&1&100\\
&(0.0008)&(0.0843)&(0)&(0)&\\
T-$\delta=0.1$&0.0023&0.2330&1&1&100\\
&(0.0008)&(0.0843)&(0)&(0)&\\\hline
\end{tabular*}}}\vspace*{10pt}
%
%
\caption{Results for the model with censored matrix, $s=2$}
\label{table2}
{\fontsize{8.8}{10.65}\selectfont{\begin{tabular*}{\tablewidth}{@{\extracolsep{\fill}}lcd{2.5}ccc@{}}
\hline
& \multicolumn{1}{c}{$\bolds{\operatorname{Err}_1}$}
& \multicolumn{1}{c}{$\bolds{\operatorname{Err}_2}$}
& \multicolumn{1}{c}{$\bolds{\operatorname{Nb}_1}$}
& \multicolumn{1}{c}{$\bolds{\operatorname{Nb}_2}$}
& \multicolumn{1}{c@{}}{\textbf{Exact}}\\
\hline
Lasso &0.1262&23.57&\multicolumn{1}{r}{96.47\phantom{00)}}&2&\phantom{00}0\\
&(0.0218)&(3.813)&\multicolumn{1}{r}{(1.670)\phantom{0}}&(0)&\\
T-Lasso&0.0456&4.688&\multicolumn{1}{r}{2.290\phantom{0)}}&2&\phantom{0}77\\
&(0.0194)&(2.157)&\multicolumn{1}{r}{(0.5881)}&(0)&\\
Dantzig &0.0792&8.000&\multicolumn{1}{r}{68.79\phantom{00)}}&2&\phantom{00}0\\
&(0.0149)&(2.159)&\multicolumn{1}{r}{(4.901)\phantom{0}}&(0)&\\
T-Dantzig&0.0404&4.0558&\multicolumn{1}{r}{2.04\phantom{00)}}&2&\phantom{0}97\\
&(0.0143)&(1.612)&(0.2416)&(0)&\\
$\delta=0.05$ &0.0064&0.6654&\multicolumn{1}{r}{2.15\phantom{00)}}&2&\phantom{0}89\\
&(0.0039)&(0.4247)&(0.4769)&(0)&\\
T-$\delta=0.05$&0.0063&0.6535&2&2&100\\
&(0.0039)&(0.4314)&(0)&(0)&\\
$\delta=0.075$ &0.0015&0.1535&2&2&100\\
&(0.0016)&(0.1637)&(0)&(0)&\\
T-$\delta=0.075$&0.0015&0.1535&2&2&100\\
&(0.0016)&(0.1637)&(0)&(0)&\\
$\delta=0.1$&0.0059&0.5410&2&2&100\\
&(0.0045)&(0.3773)&(0)&(0)&\\
T-$\delta=0.1$&0.0059&0.5410&2&2&100\\
&(0.0045)&(0.3773)&(0)&(0)&\\\hline
\end{tabular*}}}
\end{table}

\begin{table}
\caption{Results for the model with censored matrix, $s=3$}
\label{table3}
{\fontsize{8.8}{10.65}\selectfont{\begin{tabular*}{\tablewidth}{@{\extracolsep{\fill}}lcd{2.5}d{2.5}cc@{}}
\hline
& \multicolumn{1}{c}{$\bolds{\operatorname{Err}_1}$}
& \multicolumn{1}{c}{$\bolds{\operatorname{Err}_2}$}
& \multicolumn{1}{c}{$\bolds{\operatorname{Nb}_1}$}
& \multicolumn{1}{c}{$\bolds{\operatorname{Nb}_2}$}
& \multicolumn{1}{c@{}}{\textbf{Exact}}\\
\hline
Lasso &0.1834&34.54&96.91&3&\phantom{00}0\\
&(0.0326)&(6.156)&(1.407)&(0)&\\
T-Lasso&0.0776&8.832&4.28&3&\phantom{0}25\\
&(0.0306)&(3.907)&(1.068)&(0)&\\
Dantzig &0.1209&12.27&73.83&3&\phantom{00}0\\
&(0.0259)&(3.556)&(3.945)&(0)&\\
T-Dantzig&0.0597&6.108&3.40&3&\phantom{0}66\\
&(0.0251)&(2.877)&(0.6164)&(0)&\\
$\delta=0.05$ &0.0055&0.5287&3.19&3&\phantom{0}85\\
&(0.0059)&(0.4952)&(0.5038)&(0)&\\
T-$\delta=0.05$&0.0053&0.5209&\multicolumn{1}{c}{3}&3&100\\
&(0.0058)&(0.5064)&\multicolumn{1}{c}{(0)}&(0)&\\
$\delta=0.075$ &0.0148&1.296&3.05&3&\phantom{0}95\\
&(0.0110)&(0.7843)&(0.2179)&(0)&\\
T-$\delta=0.075$&0.0148&1.302&\multicolumn{1}{c}{3}&3&100\\
&(0.0109)&(0.7935)&\multicolumn{1}{c}{(0)}&(0)&\\
$\delta=0.1$&0.0415&3.791&3.02&3&\phantom{0}98\\
&(0.0177)&(1.1552)&(0.1400)&(0)&\\
T-$\delta=0.1$&0.0415&3.793&\multicolumn{1}{c}{3}&3&100\\
&(0.0177)&(1.159)&\multicolumn{1}{c}{(0)}&(0)&\\\hline
\end{tabular*}}}\vspace*{10pt}
%
%
\caption{Results for the model with censored matrix, $s=5$}
\label{table4}
{\fontsize{8.8}{10.65}\selectfont{\begin{tabular*}{\tablewidth}{@{\extracolsep{\fill}}lcd{2.4}d{2.5}cc@{}}
\hline
& \multicolumn{1}{c}{$\bolds{\operatorname{Err}_1}$}
& \multicolumn{1}{c}{$\bolds{\operatorname{Err}_2}$}
& \multicolumn{1}{c}{$\bolds{\operatorname{Nb}_1}$}
& \multicolumn{1}{c}{$\bolds{\operatorname{Nb}_2}$}
& \multicolumn{1}{c@{}}{\textbf{Exact}}\\
\hline
Lasso &0.3183&57.68&97.57&5&\phantom{0}0\\
&(0.0596)&(10.51)&(1.089)&(0)&\\
T-Lasso&0.1693&20.22&10.31&5&\phantom{0}0\\
&(0.0551)&(7.408)&(2.331)&(0)&\\
Dantzig &0.2225&22.68&81.04&5&\phantom{0}0\\
&(0.0429)&(6.275)&(3.967)&(0)&\\
T-Dantzig&0.1159&12.08&7.87&5&\phantom{0}3\\
&(0.0430)&(5.174)&(1.6891)&(0)&\\
$\delta=0.05$ &0.0596&4.544&5.52&5&63\\
&(0.0417)&(2.457)&(0.8423)&(0)&\\
T-$\delta=0.05$&0.0592&4.613&5.08&5&92\\
&(0.0414)&(2.535)&(0.2712)&(0)&\\
$\delta=0.075$ &0.1327&11.11&5.12&5&91\\
&(0.0566)&(3.059)&(0.4069)&(0)&\\
T-$\delta=0.075$&0.1327&11.14&5.03&5&97\\
&(0.0565)&(3.097)&(0.1705)&(0)&\\
$\delta=0.1$&0.2331&20.29&5.06&5&95\\
&(0.0698)&(3.154)&(0.2764)&(0)&\\
T-$\delta=0.1$&0.2371&20.61&4.97&4.95&98\\
&(0.0792)&(3.933)&(0.2628)&(0.21)&\\\hline
\end{tabular*}}}
\end{table}

\begin{table}
\caption{Results for the model with censored matrix, $s=10$}
\label{table5}
\begin{tabular*}{\tablewidth}{@{\extracolsep{\fill}}lcd{3.4}d{2.5}d{2.5}c@{}}
\hline
& \multicolumn{1}{c}{$\bolds{\operatorname{Err}_1}$}
& \multicolumn{1}{c}{$\bolds{\operatorname{Err}_2}$}
& \multicolumn{1}{c}{$\bolds{\operatorname{Nb}_1}$}
& \multicolumn{1}{c}{$\bolds{\operatorname{Nb}_2}$}
& \multicolumn{1}{c@{}}{\textbf{Exact}}\\
\hline
Lasso &0.7181&100.7&97.98&\multicolumn{1}{c}{10}&\phantom{0}0\\
&(0.1426)&(19.38)&(0.8364)&\multicolumn{1}{c}{(0)}&\\
T-Lasso&0.5560&55.09&27.02&\multicolumn{1}{c}{10}&\phantom{0}0\\
&(0.1499)&(14.57)&(3.781)&\multicolumn{1}{c}{(0)}&\\
Dantzig &0.5625&55.11&87.71&\multicolumn{1}{c}{10}&\phantom{0}0\\
&(0.1383)&(13.33)&(3.672)&\multicolumn{1}{c}{(0)}&\\
T-Dantzig&0.4203&40.41&22.33&9.98&\phantom{0}0\\
&(0.1467)&(12.91)&(3.212)&(0.1400)&\\
$\delta=0.01$ &0.3142&24.91&31.6&\multicolumn{1}{c}{10}&\phantom{0}0\\
&(0.1614)&(7.068)&(4.079)&\multicolumn{1}{c}{(0)}&\\
T-$\delta=0.01$&0.2760&20.41&14.13&9.95&\phantom{0}0\\
&(0.1612)&(7.539)&(1.677)&(0.2598)&\\
$\delta=0.05$ &0.9679&56.18&14.11&9.33&\phantom{0}2\\
&(0.3688)&(14.65)&(2.403)&(0.8724)&\\
T-$\delta=0.05$&1.0187&62.38&10.07&8.23&16\\
&(0.4088)&(17.90)&(1.226)&(1.535)&\\
$\delta=0.1$&1.392&98.89&10.31&7.94&14\\
&(0.2821)&(11.27)&(1.514)&(1.391)&\\
T-$\delta=0.1$&1.483&108.1&6.92&5.95&37\\
&(0.3003)&(12.99)&(1.324)&(1.519)&\\\hline
\end{tabular*}
\end{table}

Our first observation is that using the Lasso estimator or
the Dantzig selector (i.e., ignoring the matrix uncertainty) has
severe consequences. These methods exhibit erratic behavior already
for the minimal sparsity $s=1$. Though their sets of nonzero
components steadily include the relevant set, they are much too
large and the results are very far from the correct selection. We
also see that the MU-selector strictly improves upon the Lasso
estimator and the Dantzig selector for all the considered error
criteria and values of $s$. In particular, for $\delta=0.1$ and
$s=1,2,3,5$, it almost systematically retrieves the sparsity pattern
and the two error measures remain very small. This is obviously no
longer the case for the bigger value $s=10$. However, note that the
MU-selector remains quite satisfactory in terms of selecting the
sparsity pattern since the average number of retrieved coefficients
is about $10$ and the average number of retrieved coefficients is
about $8$. Thresholding the coefficients logically improves the
retrieved sparsity patterns of the Lasso estimator and Dantzig
selector. Nevertheless, in most of the cases the MU-selector
outperforms their thresholded versions as well. This fact is even
more significant because we simulate with a threshold which has been
well chosen knowing the true value of the nonzero coefficients. In
practice, choosing a relevant threshold is a very intricate question
since the order of magnitude of the nonzero coefficients is
typically unknown. On the other hand, for the MU-selector
thresholding can be avoided. Indeed, its effect is not significant,
especially when $s$ is small. This is due to the fact that the
original (nonthresholded) MU-selector is already very accurate in
recovering the sparsity pattern.

Finally, note that the good results for the MU-selector are not
due to the fact that we optimize over $\Theta=\mathbb{R}_{+}^{500}$
instead of $\Theta=\mathbb{R}^{500}$. In particular, taking
$\delta=0$ leads to the same kind of results as those for the
Dantzig selector.

\subsection{Model with missing data}
We consider now the model with missing data as defined in Section \ref{sec2}.
%
\begin{table}[b]
\caption{Results for the model with missing data, $s=1$}
\label{table6}
{{\begin{tabular*}{\tablewidth}{@{\extracolsep{\fill}}lcd{2.5}d{2.4}cc@{}}
\hline
& \multicolumn{1}{c}{$\bolds{\operatorname{Err}_1}$}
& \multicolumn{1}{c}{$\bolds{\operatorname{Err}_2}$}
& \multicolumn{1}{c}{$\bolds{\operatorname{Nb}_1}$}
& \multicolumn{1}{c}{$\bolds{\operatorname{Nb}_2}$}
& \multicolumn{1}{c@{}}{\textbf{Exact}}\\
\hline
Lasso &0.0212&2.606&94.59&1&\phantom{00}0\\
&(0.0105)&(1.232)&(3.256)&\multicolumn{1}{c}{(0)}&\\
T-Lasso&0.0011&0.111&\multicolumn{1}{c}{1}&1&100\\
&(0.0010)&(0.1019)&\multicolumn{1}{c}{(0)}&\multicolumn{1}{c}{(0)}&\\
Dantzig &0.0109&1.114&64.24&\multicolumn{1}{c}{1}&\phantom{00}0\\
&(0.0072)&(0.7360)&(11.06)&\multicolumn{1}{c}{(0)}&\\
T-Dantzig&0.0011&0.1097&\multicolumn{1}{c}{1}&1&100\\
&(0.0010)&(0.1030)&\multicolumn{1}{c}{(0)}&\multicolumn{1}{c}{(0)}&\\
$\delta=0.05$ &0.0041&0.3376&8.55&\multicolumn{1}{c}{1}&\phantom{00}6\\
&(0.0029)&(0.2218)&(5.087)&\multicolumn{1}{c}{(0)}&\\
T-$\delta=0.05$&0.0022&0.2271&\multicolumn{1}{c}{1}&1&100\\
&(0.0012)&(0.1203)&\multicolumn{1}{c}{(0)}&\multicolumn{1}{c}{(0)}&\\
$\delta=0.075$ &0.0039&0.3449&3.99&\multicolumn{1}{c}{1}&\phantom{0}29\\
&(0.0021)&(0.1625)&(3.090)&\multicolumn{1}{c}{(0)}&\\
T-$\delta=0.075$&0.0031&0.3133&\multicolumn{1}{c}{1}&1&100\\
&(0.0011)&(0.1124)&\multicolumn{1}{c}{(0)}&\multicolumn{1}{c}{(0)}&\\
$\delta=0.1$&0.0047&0.4490&1.94&\multicolumn{1}{c}{1}&\phantom{0}61\\
&(0.0019)&(0.1356)&(1.605)&\multicolumn{1}{c}{(0)}&\\
T-$\delta=0.1$&0.0044&0.4451&\multicolumn{1}{c}{1}&1&100\\
&(0.0012)&(0.1268)&\multicolumn{1}{c}{(0)}&\multicolumn{1}{c}{(0)}&\\\hline
\end{tabular*}}}
\end{table}
\begin{table}
\caption{Results for the model with missing data, $s=2$}
\label{table7}
{\fontsize{8.7}{10.6}\selectfont{\begin{tabular*}{\tablewidth}{@{\extracolsep{\fill}}lcd{2.5}d{2.5}cd{3.0}@{}}
\hline
& \multicolumn{1}{c}{$\bolds{\operatorname{Err}_1}$}
& \multicolumn{1}{c}{$\bolds{\operatorname{Err}_2}$}
& \multicolumn{1}{c}{$\bolds{\operatorname{Nb}_1}$}
& \multicolumn{1}{c}{$\bolds{\operatorname{Nb}_2}$}
& \multicolumn{1}{c@{}}{\textbf{Exact}}\\
\hline
Lasso &0.0425&4.675&96.02&2&0\\
&(0.0162)&(1.786)&(2.074)&\multicolumn{1}{c}{(0)}&\\
T-Lasso&0.0047&0.4481&2.02&\multicolumn{1}{c}{2}&98\\
&(0.0037)&(0.3318)&(0.1400)&\multicolumn{1}{c}{(0)}&\\
Dantzig &0.0269&2.695&74.39&\multicolumn{1}{c}{2}&0\\
&(0.0134)&(1.4823)&(5.774)&\multicolumn{1}{c}{(0)}&\\
T-Dantzig&0.0046&0.4330&\multicolumn{1}{c}{2}&2&100\\
&(0.0035)&(0.3194)&\multicolumn{1}{c}{(0)}&\multicolumn{1}{c}{(0)}&\\
$\delta=0.05$ &0.0131&1.033&6.76&\multicolumn{1}{c}{2}&13\\
&(0.0078)&(0.4688)&(3.572)&\multicolumn{1}{c}{(0)}&\\
T-$\delta=0.05$&0.0106&1.018&\multicolumn{1}{c}{2}&2&100\\
&(0.0055)&(0.4692)&\multicolumn{1}{c}{(0)}&\multicolumn{1}{c}{(0)}&\\
$\delta=0.075$ &0.0167&1.517&3.20&\multicolumn{1}{c}{2}&48\\
&(0.0071)&(0.4557)&(1.489)&\multicolumn{1}{c}{(0)}&\\
T-$\delta=0.075$&0.0160&1.525&2.01&\multicolumn{1}{c}{2}&99\\
&(0.0064)&(0.4584)&(0.099)&\multicolumn{1}{c}{(0)}&\\
$\delta=0.1$&0.0247&2.351&2.27&\multicolumn{1}{c}{2}&77\\
&(0.0074)&(0.4634)&(0.5264)&\multicolumn{1}{c}{(0)}&\\
T-$\delta=0.1$&0.0245&2.362&\multicolumn{1}{c}{2}&2&100\\
&(0.0070)&(0.4731)&\multicolumn{1}{c}{(0)}&\multicolumn{1}{c}{(0)}&\\\hline
\end{tabular*}}}\vspace*{9pt}
%
\caption{Results for the model with missing data, $s=3$}
\label{table8}
{\fontsize{8.7}{10.60}\selectfont{\begin{tabular*}{\tablewidth}{@{\extracolsep{\fill}}lcd{2.5}d{2.5}cd{3.0}@{}}
\hline
& \multicolumn{1}{c}{$\bolds{\operatorname{Err}_1}$}
& \multicolumn{1}{c}{$\bolds{\operatorname{Err}_2}$}
& \multicolumn{1}{c}{$\bolds{\operatorname{Nb}_1}$}
& \multicolumn{1}{c}{$\bolds{\operatorname{Nb}_2}$}
& \multicolumn{1}{c@{}}{\textbf{Exact}}\\
\hline
Lasso &0.0721&6.828&96.89&3&0\\
&(0.0251)&(2.116)&(1.449)&(0)&\\
T-Lasso&0.0134&1.225&3.12&3&88\\
&(0.0093)&(0.8126)&(0.3250)&(0)&\\
Dantzig &0.0496&4.844&80.45&3&0\\
&(0.0204)&(2.075)&(4.693)&(0)&\\
T-Dantzig&0.0117&1.119&3.05&3&95\\
&(0.0082)&(0.8438)&(0.2180)&(0)&\\
$\delta=0.05$ &0.0322&2.591&6.8&3&10\\
&(0.0138)&(0.7730)&(2.942)&(0)&\\
T-$\delta=0.05$&0.0293&2.726&3.04&3&96\\
&(0.0119)&(0.8735)&(0.1959)&(0)&\\
$\delta=0.075$ &0.0439&4.0308&3.96&3&50\\
&(0.0137)&(0.7988)&(1.333)&(0)&\\
T-$\delta=0.075$&0.0432&4.1098&3.01&3&99\\
&(0.0130)&(0.8505)&(0.0994)&(0)&\\
$\delta=0.1$&0.0653&6.217&3.21&3&84\\
&(0.0160)&(0.8355)&(0.5156)&(0)&\\
T-$\delta=0.1$&0.0651&6.235&\multicolumn{1}{c}{3}&3&100\\
&(0.0158)&(0.8500)&\multicolumn{1}{c}{(0)}&(0)&\\\hline
\end{tabular*}}}
\end{table}
We design the numerical experiment in the same way as in Section
\ref{cens} except that the observed matrix $Z$ is now given by
(\ref{ex1}) with $\pi=0.1$.

\subsubsection*{Results}

The results are given in Tables \ref{table6}--\ref{table10}.
We see that again the Lasso and Dantzig selector are
highly unstable in selecting the sparsity pattern, whereas the
MU-selector does a good job. The thresholded estimators T-Lasso
and T-Dantzig are also quite accurate in retrieving the sparsity
pattern, except for $s=10$. However, in all the cases the
MU-selector does it better. The MU-selector with $\delta=0.05$
(or $\delta=0.01$ for $s=10$) has the smallest error measures
$\operatorname{Err}_1$ and $\operatorname{Err}_2$, whereas the
sparsity pattern is
%
\begin{table}
\caption{Results for the model with missing data, $s=5$}
\label{table9}
{\fontsize{8.8}{10.65}\selectfont{\begin{tabular*}{\tablewidth}{@{\extracolsep{\fill}}lcd{2.4}d{2.5}cd{2.0}@{}}
\hline
& \multicolumn{1}{c}{$\bolds{\operatorname{Err}_1}$}
& \multicolumn{1}{c}{$\bolds{\operatorname{Err}_2}$}
& \multicolumn{1}{c}{$\bolds{\operatorname{Nb}_1}$}
& \multicolumn{1}{c}{$\bolds{\operatorname{Nb}_2}$}
& \multicolumn{1}{c@{}}{\textbf{Exact}}\\
\hline
Lasso &0.1302&9.993&97.24&5&0\\
&(0.0499)&(2.657)&(1.097)&(0)&\\
T-Lasso&0.0418&3.331&5.65&5&56\\
&(0.0326)&(2.056)&(0.899)&(0)&\\
Dantzig &0.1005&9.371&84.36&5&0\\
&(0.0443)&(4.113)&(4.009)&(0)&\\
T-Dantzig&0.0365&3.356&5.38&5&74\\
&(0.0275)&(2.454)&(0.7454)&(0)&\\
$\delta=0.05$ &0.1033&8.301&8.19&5&14\\
&(0.0384)&(1.713)&(2.591)&(0)&\\
T-$\delta=0.05$&0.1001&8.900&5.14&5&87\\
&(0.0362)&(2.146)&(0.3746)&(0)&\\
$\delta=0.075$ &0.1485&13.25&5.96&5&48\\
&(0.0415)&(1.716)&(1.272)&(0)&\\
T-$\delta=0.075$&0.1477&13.53&5.05&5&95\\
&(0.0402)&(1.999)&(0.2179)&(0)&\\
$\delta=0.1$&0.2133&19.60&5.31&5&79\\
&(0.0494)&(1.708)&(0.7835)&(0)&\\
T-$\delta=0.1$&0.2131&19.70&5.03&5&97\\
&(0.0488)&(1.904)&(0.1705)&(0)&\\\hline
\end{tabular*}}}\vspace*{10pt}
%
\caption{Results for the model with missing data, $s=10$}
\label{table10}
{\fontsize{8.8}{10.65}\selectfont{\begin{tabular*}{\tablewidth}{@{\extracolsep{\fill}}ld{2.5}d{2.4}d{2.5}d{2.5}d{2.0}@{}}
\hline
& \multicolumn{1}{c}{$\bolds{\operatorname{Err}_1}$}
& \multicolumn{1}{c}{$\bolds{\operatorname{Err}_2}$}
& \multicolumn{1}{c}{$\bolds{\operatorname{Nb}_1}$}
& \multicolumn{1}{c}{$\bolds{\operatorname{Nb}_2}$}
& \multicolumn{1}{c@{}}{\textbf{Exact}}\\
\hline
Lasso &0.4746&18.04&98.01&\multicolumn{1}{c}{10}&0\\
&(0.1702)&(4.334)&(0.7549)&\multicolumn{1}{c}{(0)}&\\
T-Lasso&(0.1710)&(3.357)&(5.317)&\multicolumn{1}{c}{(0)}&0\\
&0.4358&13.74&37.44&\multicolumn{1}{c}{10}&\\
Dantzig &0.4229&38.13&90.77&\multicolumn{1}{c}{10}&0\\
&(0.1684)&(15.95)&(2.853)&\multicolumn{1}{c}{(0)}&\\
T-Dantzig&0.3862&34.84&32.77&\multicolumn{1}{c}{10}&0\\
&(0.1690)&(16.35)&(5.184)&\multicolumn{1}{c}{(0)}&\\
$\delta=0.01$ &0.2891&10.78&47.38&\multicolumn{1}{c}{10}&0\\
&(0.1285)&(2.059)&(5.351)&\multicolumn{1}{c}{(0)}&\\
T-$\delta=0.01$&0.2725&12.26&19.93&\multicolumn{1}{c}{10}&0\\
&(0.1271)&(2.719)&(3.311)&\multicolumn{1}{c}{(0)}&\\
$\delta=0.05$ &0.7710&45.89&18.02&9.91&0\\
&(0.2755)&(6.212)&(3.720)&(0.2861)&\\
T-$\delta=0.05$&0.7719&48.36&13.41&9.73&6\\
&(0.2807)&(7.134)&(2.015)&(0.6611)&\\
$\delta=0.1$&1.182&84.80&13.42&9.37&6\\
&(0.2983)&(8.477)&(2.324)&(0.8204)&\\
T-$\delta=0.1$&1.196&87.42&10.81&8.78&23\\
&(0.3104)&(9.304)&(1.521)&(1.338)&\\\hline
\end{tabular*}}}
\end{table}
better retrieved for $\delta=0.1$. This reflects a tradeoff between
estimation and selection. Smaller values of $\delta$ lead to smaller
errors $\operatorname{Err}_1$ and $\operatorname{Err}_2$, whereas
larger values of
$\delta$ lead to a very accurate recovery of the sparsity pattern.
The error measures $\operatorname{Err}_1$ and $\operatorname{Err}_2$
of the
thresholded estimators T-Lasso and \mbox{T-Dantzig} are somewhat smaller
than those of the MU-selector, except for $s=10$. Note, however,
that we report the results for the performance of T-Lasso and
T-Dantzig with a threshold based on the knowledge of the true
coefficients.

\subsection{Portfolio replication}

We now present a ``toy'' application based on financial data. We
apply model (\ref{0}) and (\ref{00}) and the MU-selector in the
context of portfolio replication as described in Section \ref{sec2}. We take
the data of the open and close prices of $p=491$ assets in the
Standard and Poors S\&P 500 index for the $n=251$ trading days of
2007. These data are provided by the Yahoo Finance Database. The
assets we use are those available for the whole year.

\subsubsection*{Experiment} Let $p^o_{ij}$ and $p^c_{ij}$ denote the open
and close prices of the $j$th asset for the $i$th day. Our experiment
is the following.\vspace*{1pt}

-- We consider the matrix $\tilde{X}$ with entries
$(\tilde{X})_{ij}=p^c_{ij}-p^o_{ij}$ and define $X$ as the
normalized matrix obtained from $\tilde{X}$.

-- We pick $s$ assets to build our portfolio. The
coordinate of each chosen asset in the vector
$\theta\in\mathbb{R}^{491}$ is set to $1/s$ and the other
coordinates to 0 [note that, in practice, if the $j$th asset is in
the portfolio, it means that the corresponding coordinate of $\theta$ is
$1/(s\tilde{\sigma_j})$, where $\tilde{\sigma_j}$ is the empirical
standard deviation of its absolute returns].

-- We consider six portfolios (see Table \ref{table11}).

%
\begin{table} 
\tablewidth=240pt
\caption{Initial portfolios}
\label{table11}
\begin{tabular*}{\tablewidth}{@{\extracolsep{\fill}}ll@{}}
\hline
$\bolds{s=2}$ & \multicolumn{1}{c@{}}{$\bolds{s=3}$} \\
\hline
Boeing, Goldman Sachs & Boeing, Google, Goldman Sachs\\
Boeing, Coca Cola & Boeing, Google, Coca Cola \\
Boeing, Ford & Boeing, Google, Ford \\
\hline
\end{tabular*}
\end{table}

-- We compute $y=X\theta+\xi$ where $\xi$ is the same
noise as in Section \ref{cens}. In practice, the noise $\xi$ can
reflect
an uncertainty about the management costs, a lack of transparency in
the definition of the returns of the portfolio or some rounding approximations.

-- We consider a matrix uncertainty of the following
type: $Z$ is obtained from $X$ by replacing one of its columns by
the zero column. The column corresponds to one of the assets in the
portfolio. The goal of this manipulation is to mimic the fact that
in practice not all the existing assets are in our restricted class.
One of the assets in the portfolio does not belong to the restricted
class since the corresponding column of $X$ is suppressed. Of
course, this asset cannot be retrieved. We suppress the column
associated to an asset different from
Boeing and Google.

-- We solve (\ref{9}) with such a matrix $Z$, with
$\delta=0.5$ and $\varepsilon$ chosen as in Section~\ref{cens}. We also
compute the Lasso estimator and the Dantzig selector.

\subsubsection*{Results}
We write B for Boeing and G for Google. The initial portfolios and
the portfolios retrieved by the MU-selector are presented in
Table \ref{table12}.

\begin{table}
\tablewidth=240pt
\caption{Retrieved portfolios, MU-selector}
\label{table12}
\begin{tabular*}{\tablewidth}{@{\extracolsep{\fill}}ll@{}}
\hline
\textbf{Initial portfolio} & \multicolumn{1}{c@{}}{\textbf{Retrieved portfolio}}\\
\hline
B, Goldman Sachs & B, Morgan Stanley, Merrill Lynch\\
B, Coca Cola & B, Pepsico\\
B, Ford & B, General Motors\\
B, G, Goldman Sachs & B, G, Morgan Stanley, Merrill Lynch\\
B, G, Coca Cola & B, G\\
B, G, Ford & B, G, General Motors\\
\hline
\end{tabular*}
\end{table}

The results are very satisfying. Indeed, the algorithm
almost always finds the correct number of assets in the portfolio
and the discarded asset is replaced by one or two assets that are
intuitively close to it. Moreover, if one takes $\delta=0.4$, then for
the initial portfolio (Boeing, Google, Coca Cola) the retrieved
portfolio becomes (Boeing, Google, Pepsico), whereas the other
results remain the same. Finally, note that the Lasso estimator and
the Dantzig selector (usual Dantzig selector or MU-selector with
$\delta=0$) systematically output more than 20 assets in the
retrieved portfolio.

\printaddresses

\end{document}